\documentclass[10pt,leqno]{amsart}
\usepackage{amssymb, amstext, amsthm, amscd, amsmath,color}
\usepackage{graphicx} 
\usepackage{tikz-cd}
\usepackage{cancel}
\usepackage{thmtools}
\usepackage{stmaryrd}
\usepackage{comment}
\usepackage{enumitem}
\usetikzlibrary{decorations.pathmorphing,shapes}
\usepackage{hyperref}
\usepackage{float}

\usepackage[textwidth=0.89in,textsize=scriptsize]{todonotes}

\usepackage{indentfirst,csquotes}

\hypersetup{ colorlinks=true, linkcolor=black, filecolor=black, urlcolor=black }

\usepackage{lipsum}

\hypersetup{
colorlinks=false,
linktoc=all
}

\theoremstyle{plain}
\newtheorem{thm}{Theorem}[subsection]%Section if you want to number it by section]
\newtheorem{cor}[thm]{Corollary}
\newtheorem{prop}[thm]{Proposition}
\newtheorem{lem}[thm]{Lemma}

\theoremstyle{definition}
\newtheorem{rem}[thm]{Remark}

\newtheorem{rems}[thm]{Remarks}
\newtheorem{defn}[thm]{Definition}

\newtheorem{eg}[thm]{Example}
\newtheorem{egs}[thm]{Examples}

\newcommand{\fty}{\infty}

\newcommand{\mR}{{\mathbb R}}
\newcommand{\mZ}{{\mathbb Z}}
\newcommand{\mN}{{\mathbb N}}

\newcommand{\ZB}{{\mathcal Z \mathcal B}}
\newcommand{\cZ}{{\mathcal Z}}
\newcommand{\cA}{{\mathcal A}}
\newcommand{\cB}{{\mathcal B}}
\newcommand{\cP}{{\mathcal P}}
\newcommand{\cC}{{\mathcal C}}

\newcommand{\id}{\textbf{id}}

\newcommand{\Cat}{\mathbf{Cat}}

\newcommand{\mbA}{\mathbf{A}}

\newcommand{\mbG}{\mathbf{G}}

\newcommand{\ChA}{{\mathbf{Ch}(\cA)}}

\newcommand{\Set}{{\mathbf {Set}}}

\newcommand{\PoSet}{{\mathbf{PoSet}}}

\newcommand{\Sub}{{\mathbf {Sub}}}

\newcommand{\zb}{{\mathbf {zb}}}
\newcommand{\mbX}{{\mathbf {X}}}

\newcommand{\noi}{{\noindent}}

\newsavebox{\pullback}
\sbox\pullback{%
\begin{tikzpicture}%
\draw (0,0) -- (1ex,0ex);%
\draw (1ex,0ex) -- (1ex,1ex);%
\end{tikzpicture}}
\newsavebox{\pushout}
\sbox\pushout{%
\begin{tikzpicture}%
\draw (0,0) -- (0ex,1ex);%
\draw (0ex,1ex) -- (1ex,1ex);%
\end{tikzpicture}}
\newsavebox{\urpushout}
\sbox\urpushout{%
\begin{tikzpicture}%
\draw (0,0) -- (1ex,0ex);%
\draw (0ex,0ex) -- (0ex,1ex);%
\end{tikzpicture}}
\newsavebox{\vpullback}
\sbox\vpullback{%
\begin{tikzpicture}%
\draw (0ex,1ex) -- (1ex, 0ex);%
\draw (1ex, 0ex) -- (2ex,1ex);%
\end{tikzpicture}}
\newsavebox{\vpushout}
\sbox\vpushout{%
\begin{tikzpicture}%
\draw (0ex,1ex) -- (1ex, 1ex);%
\draw (1ex, 1ex) -- (2ex,1ex);%
\end{tikzpicture}}
\newcounter{sarrow}

\usepackage{xcolor}
\usepackage{soul}

\begin{document}

\title{Change Action Derivatives in Persistent Homology} %%%%%%%%%%%%
\author[DS]{Deni Salja}
\date{\today}
\address{Address}
\email{denisalja@dal.ca}
\maketitle

\section*{Introduction}

Persistent homology is a popular technique in topological data analysis that tracks the lifespans of homological features in a nested sequence of spaces. This data is typically presented in a multi-set called a persistence diagram or a barcode \cite{CZ-ComputingMultiPersistence}, \cite{peristenceBarcodesforShapes}, and \cite{E&H}.
When a topological space is filtered with respect to a single parameter in $\mR$, and when homology coefficients are taken in a principal ideal domain (PID), the persistence diagram can be computed using the presentation theorem for finitely generated modules over a PID \cite{CZ-ComputingMultiPersistence}. The modules corresponding to multi-parameter filtrations, whose indexing poset is $\mR^n$ for some $n \geq 2$, lack a similar presentation theorem which means no complete invariant for these exists \cite{CZ-ComputingMultiPersistence}. 

One way to extract the number of homology classes with a given interval lifespan is to calculate the rank of the pair-group as mentioned in \cite[Section 2 - Tame Functions]{E&H}.
The pair group at a given time interval is a quotient of homology classes whose lifespans have already passed by the classes whose lifespans passed sometime before.
The rank of this quotient counts the number of homology classes as a difference in an abelian group. By scanning through the domain of the filtration and calculating the ranks of pair groups, one can recover the data of a persistence diagram/barcode.
The purpose of this paper is to connect persistent homology to a categorical notion of differentiation which provides semantics for incremental computation.

\section*{Notation and Structure} 

Categories are generally written in bold font: $\mathbf{A}$, $\mathbf{B}$, $\mathbf{C}$ \dots and their objects are written in capitals: $A, B,$ $C$\dots 
Juxtaposition of morphisms denotes diagrammatic composition, so for a pair of morphisms $f$ and $g$, by $fg$ we mean `$f$ first, followed by $g$.'
If applicative notation is ever used it will be explicitly mentioned or denoted with a $\circ$. 
We use $!$ to denote the unique map into a terminal object. 
The groups in this paper are generally abelian and $\mathbf{B}G$ denotes the delooping of a group $G$. 
The arrow category of a category $\cC$ is denoted $\mathbf{Arr}(\cC)$. 

The first two sections of the paper introduce change-actions, change-action derivatives, and finite differences with relevant examples. 
The main change-action derivative of interest in this paper is a functor and recognizing it requires a categorification of the calculus of finite differences from \cite[Section 5.2]{AP&L}, said categorification appears in Theorem \ref{thm functors into G1 are differentiable}. 

A functorial framework for persistent homology is described in the third section, where we define filtrations as functors on a finite poset valued in chain complexes of an abelian group.
Comparing cycles and boundaries that appear at pairs of comparable points in a poset can be impractical and full of redundant information so we extend the cycle and boundary functors to the topology of upward closed subsets.
Regions in the indexing poset are then represented by nested pairs of open subsets which are given by objects in the arrow category of (the opposite category of) the open lattice.
The appearance and disappearance of homology data in the filtration is encoded by the `homological memory' and `homological lifespan' functors, $\ZB Fn$ and $\Gamma_n$, which are defined as functors on the aforementioned arrow category in Sections \ref{SS homological memory} and \ref{SS Lifespan subquot} respectively. 

A change-action structure on the opens in the topology is given by `shifting blankets' in Section \ref{SS blanket shift change action} and the homological memory functor $\ZB F_n$ is extended to fit together with this change-action structure in Section \ref{S CAD in multipersistence}.
Viewing the rank function for an abelian category as a functor on the arrow category, we compose it with $\Gamma_n$ to get a generalization of the rank of the pair group from \cite[Section 2 - Tame Functions]{E&H}.
The categorified calculus of finite differences applies to the rank $\ZB F_n$ in the main result, Theorem \ref{thm extended homological memory functor is differentiable} and Corollary \ref{cor lifespan functor recovered by derivative}, giving a change-action derivative that recovers the rank of $\Gamma_n$. 

\section*{Acknowledgements} 

This project was the topic of a summer research project and my honours thesis at UCalgary, supervised by Kristine Bauer. 
It was funded by UCalgary's Program for Undergraduate Research Experience with portions of it stemming from a previous summer project funded by NSERC's USRA program.
A big thank you goes out to Kristine for her mentorship during my time at UCalgary and as well as Geoff Vooys for helpful discussions, feedback, and encouraging me to get this off my laptop and into the world. 

%Future work on this topic involves exploring different topological invariants with more general coefficients, namely persistent sheaf cohomology, and a relationship to functor calculus. 

%SECTION 1

\section{Change Actions and Derivatives}

Change structures were proposed in \cite{Cai2013ATO} to give a semantic framework for incremental differentiation.
The change-actions introduced by \cite{A-P&0)} extend these ideas to monoidal categories and describe a categorical notion of (higher order) generalized differentiation that we use to describe the incremental changes in persistent homology. 

Change actions are defined in terms of monoid actions internal to a category, and hence require the definition of a monoid object.
While monoids and monoid actions can be defined in the general context of a monoidal category \cite[VII]{CFTWM}, for this paper it suffices to consider categories in which the monoidal structure is strict and given by products.
Such categories are called (strict) cartesian monoidal categories.
We use $\times$ to denote products and $\top$ to denote terminal objects (the empty products). 

\begin{egs}[(Strict) Cartesian (Monoidal) Categories] \label{CartesianMonoidalCats}\
\begin{itemize}
\item The trivial category with one objects and a single identity morphism, $\mathbf{1}$.
\item The category of sets and functions, $\Set$. % has finite products given by cartesian products and any singleton set is a terminal object.
\item The category of categories and functors, $\Cat$.
\item The category of partially ordered sets and order preserving functions, $\PoSet$.
\end{itemize}
\end{egs}

\noi An action of a monoid, $M$,  on an object, $A$,  in $\mbX$ is a map, $\alpha : A \times M \to A$ satisfying an associativity and unit condition.
The object $M$ can be thought of as representing possible changes that the action $\alpha$ can apply to $A$.
Asking for these changes to be composablable leads to the definition of a monoid. 

\begin{defn}[{\cite[III.6, page 75]{CFTWM}}]\label{Def monoid}
In a category, $\mbX$, with finite products and a terminal object, a \textit{monoid} is a triple ($M , \mu, e$) such that  $M$ is an object of $\mbX$, 
\[ 
\mu : M \times M \to M, \qquad  \qquad e : \top \to M 
\]
\noi are maps in $\mbX$,  and the following unit and associativity diagrams,
\begin{equation}
\begin{tikzcd}[column sep = large, row sep = large]
M \rar["\langle 1_M {,} !_M e\rangle "] \ar[dr, equals] \dar["\langle !_M e {,} 1_M\rangle"'] & M  \times M \dar["\mu"] \\
M \times M \rar["\mu"'] & M
\end{tikzcd}\qquad \qquad
\begin{tikzcd}[column sep = large, row sep = large]
M \times M \times M \rar["1_M \times \mu"]\dar["\mu \times 1_M"'] & M \times M \dar["\mu"]\\
M \times M \rar["\mu"'] & M 	
\end{tikzcd},
\end{equation}
commute in $\mbX$.  
\end{defn}\

Now we can state Alvarez-Picallo's definition for the discrete derivative introduced earlier in this section.  

\begin{defn}[Alvarez-Picallo \cite{A-P&0)}]\label{Def change action}
A \textit{change action} in a cartesian monoidal category $\mbX$ is an action, $\oplus_A$, of a monoid, $(\Delta A , +_A , 0_A)$, on an object $A$. In other words, a change action is a tuple
\[ 
\bar{A} = (A, \Delta A, \oplus_A, +_A, 0_A), 
\]
such that $A, \Delta A $ are objects in $ \mbX$,

\begin{center}
\begin{tikzcd}
A \times \Delta A \rar["\oplus_A"] & A ,&
\Delta A \times \Delta A \rar["+_A"] & \Delta A, &
\top \rar["0_A"] & \Delta A
\end{tikzcd}, 	
\end{center}
are maps in $\mbX$, $(\Delta A, +_A, 0_A)$ is a monoid in $\mbX$, and the following diagrams
%\[ \langle 1_A, !_A 0_A \rangle \oplus_A = 1_A \quad \text{ and } \quad (1_A \times +_A) \oplus_A = (\oplus_A \times 1_{\Delta A}) \oplus_A \]\
\begin{equation} 
\begin{tikzcd} [column sep = large, row sep = large]
 A \ar[dr, equals] \rar["\langle 1_A{,} !_A 0_A \rangle"] & A \times \Delta A \dar["\oplus"] \\
 & A	
 \end{tikzcd} 
 \quad \quad \quad 
 \begin{tikzcd} [column sep = large, row sep = large]
	A \times \Delta A \times \Delta A \dar["(1_A) \times (+_A)"']\rar["(\oplus_A) \times (1_{\Delta A})"] & A \times \Delta A \dar["\oplus_A"] \\
	 A \times \Delta A \rar["\oplus_A"'] & A 
\end{tikzcd},
\end{equation}
commute. 
\end{defn}
The left diagram says the identity of the monoid acts trivially while the diagram on the right says the action, $\oplus_A$, is compatible with the monoid operation. A change action structure on an object is a (right) monoid action on that object.

\begin{eg}\label{Eg (N,N,+,+,0)}
That natural numbers, $\mN$,  with their usual ordering is an object in $\PoSet$. The quadruple $(\mN, \mN, +,+,0)$ is a change action in the category $\PoSet$. 
\end{eg}

\begin{eg}
Every monoid is an example of a change action structure in two ways. The first action is only by identities 
\[ \begin{tikzcd} 
M \times M \rar["\pi_1"] & M
\end{tikzcd} \]
and the second is by the monoid operation 

\[ \begin{tikzcd} 
M \times M \rar["\mu"] & M
\end{tikzcd} \]
\end{eg}
Change actions capture a discrete notion of differentiation. A change-action derivative for a morphism between change-action structures tracks the incremental change in the codomain corresponding to an incremental change in the domain. 

\begin{defn}[{\cite[Definition 3.3]{AP&L}}] \label{Def derivative}
Let $(A, \Delta A, \oplus_A, +_A, 0_A)$ and $(B, \Delta B, \oplus_B, +_B, 0_B)$ be change actions in a cartesian monoidal category $\mbX$. 
A \textit{change-action derivative} of a map $f: A \to B$ in $\mbX$ is a map 
\begin{center}
\begin{tikzcd}
A \times \Delta A \rar["\partial f"] & \Delta B	
\end{tikzcd}
\end{center}
such that for any three maps
\begin{center}
\begin{tikzcd}
	X \rar["x"] & A, & X \rar[shift left, "y"] \rar[shift right, "z"'] & \Delta A 
\end{tikzcd}
\end{center}
in $\mbX$, the following axioms hold. 
\begin{enumerate}
\item (\textbf{CAD1})\label{CAD1} 
The diagram, 
\[ 
\begin{tikzcd}
	X \rar["\langle x{,} y \rangle"] \dar[ "\langle x{,} \langle x{,} y \rangle \rangle"'] & A \times \Delta A \dar["\oplus_A"] \\
	A \times ( A \times \Delta A) \dar["f \times \partial f"']& A \dar["f"] \\
	B \times \Delta B \rar["\oplus_B"] & B  
 \end{tikzcd},
\]
commutes.

\item (\textbf{CAD2}) \label{CAD2} The diagrams,
\[ 
\begin{tikzcd} 
 X \dar["\langle \langle x {,} y \rangle {,} \langle \langle x{,}y \rangle \oplus_A {,} z \rangle \rangle"'] \rar["\langle x {,} \langle y{,}z \rangle +_A \rangle"] & A \times \Delta A \dar[dd, "\partial f"] \\
 ( A \times \Delta A) \times (A \times \Delta A) \dar["\partial f \times \partial f"'] & \\
\Delta B \times \Delta B \rar["+_B"'] & \Delta B 
 \end{tikzcd} \qquad , \qquad 
 \begin{tikzcd}[column sep = large, row sep = large]
 X \rar["\langle x{,} !_X 0_{\Delta A} \rangle"] \dar["!_X"'] & A \times \Delta A \dar["\partial f"] \\
 0 \rar["0_{\Delta B}"'] & \Delta B	
 \end{tikzcd}
 \]
commute.  
\end{enumerate}
\end{defn}
The first condition is an analogue for the linear approximation of a function. In applicative notation for composition, the commuting diagram in \textbf{CAD1} becomes:
%if we take $A = B = \mR$ with both change actions given by letting $\mN$ act on $\mR$ via addition,  $X = \{ * \}$, and $x, y, z, f$ to be maps in $\PoSet$,  then the equation 
\begin{align}\label{eq linear approx of f}
f(x \oplus_A y) = f(x) \oplus_B \partial f(x, y)
\end{align}
If we take $A = B = \mR$ with both change actions given addition on $\mR$ and $X = \{ * \}$, then $x, y, z$ are functions $\{*\} \to \mR$ and $f$ is a function $\mR \to \mR$. Then equation \ref{eq linear approx of f} becomes precisely 
\[ 
f(x + y) = f(x) + \partial f (x , y) 
\]
where $\partial f (x,y)$ is the change in the codomain of $f$ (at $x$) determined by the action of $y$ in the domain (at $x$). When $f$ is differentiable at $x$, the definition $\partial f (x,y) = y f'(x)$  makes sense for sufficiently small $y$:
\[ 
f(x + y) \approx f(x) + y f '(x) 
\]
In differential calculus the derivative is always linear, so one might expect the equation
\[ \partial f ( x, y +_A z ) = \partial f ( x, y ) +_B \partial f ( y, z) \]
from the left diagram of \textbf{CAD2}, but this linearity fails precisely by an action in the first component: 
\[ \partial f ( x, y +_A z ) = \partial f ( x, y ) +_B \partial f ( x \oplus_A y, z).\]
To be clear, the left diagram in \textbf{CAD2} says change-action derivatives can be thought of as \textit{linear up to an action in the first component}. 
For example, for a function $f: \mR \to \mR$, we can take $\partial f$ to be the difference $\partial f (x,y) = f(x+y) - f(x)$ and see that  
\begin{align*}
\partial f (a , b+c) 
&=  f(a+(b+c)) - f(a) \\
&= \big( f((a+b)+c)) - f(a+b) \big) + \big( f(a+b) - f(a) \big) \\
&= \partial f (a +b , c) + \partial f (a,b) 
\end{align*}
The second part of \textbf{CAD2} says that the action by the unit in the monoid induces a trivial change, i.e. no change at all.
\[ \partial f (x , 0_{\Delta A}) = 0_{\Delta B}\]

\begin{eg}\label{Eg a derivative of k-translation on N in Poset}
Let $\mN$ be equipped with the addition action from Example \ref{Eg (N,N,+,+,0)} in $\PoSet$. Fix $k \in \mN$ and let $T_k : \mN \to \mN$ be defined by $ T_k(n) = n + k$ be the order preserving map that translates by $k$. Then the projection, 
\[ 
\pi_1 : \mN \times \mN \to \mN \ ; \ (a,b) \mapsto b, 
\]
is a change-action derivative of $T_k$. For any $a,b \in \mN$,  \textbf{CAD1} is the equation
\[ T_k(a+b) = (a+b)+k = (a+k) +b = T_k(a) + \pi_1(a,b). \]
For any $a,b,c \in \mN$, \textbf{CAD2} is the pair of equations,
\[ \pi_1(a, b+c) = b+c = \pi_1(a,b) + \pi_1(a+b,c)\]
% = \pi_1 \times \pi_1\big( (a,b) , (a+b, c) \big) \]\
and 
\[ \pi_1(a, 0) = 0 ,\]
respectively. Here $\pi_1$ accepts a pair $(a,b)$ as input and it outputs $b$, the amount that $T_k$ changes between $a$ and $a+b$. 
Since $T_k$ is monotone, $\pi_1$ can be expressed as a finite difference:
\[ \pi_1(a,b) = b = ((a+b)+k) - (a+k) = T_k(a+b) - T_k(a) . \]
\end{eg}

\begin{defn}\label{Def change action diffbl}
For any two objects $A$, $B$ in a cartesian monoidal category,  with respective monoid actions,  a map   $f: A \to B$ is said to be \textit{(change-action) differentiable (in $\mbX$ with respect to given monoid actions)} if the map $f$ has a change-action derivative as defined in \ref{Def derivative}. 
\end{defn}

The map $T_k \in \PoSet_1$ from Example \ref{Eg a derivative of k-translation on N in Poset} is differentiable in $\PoSet$.
Consider a piecewise staircase function $f : \mR \to \mR$ : for example 
\[ f(x) = \begin{cases} 
n & n \leq x < n+\frac{1}{2} \\
n + \dfrac{x - (n + \frac{1}{2})}{2} & n + \frac{1}{2} \leq x < n + 1
\end{cases},\]
with the pieces smoothed out so that it is differentiable and monotone increasing.
Although its derivative is strictly non-negative, it is not monotone since for each $n \in \mZ$  the derivative of each piece is zero in the first half of each interval $(n, n+\frac{1}{2})$ and then is equal to $\frac{1}{2}$ in the second half of each interval $(n + \frac{1}{2} , n + 1)$.
A change-action derivative of a map needs to lie in the same category as the original map.
In particular, change-action derivatives of monotone functions need to be monotone themselves, and this leads to certain monotone functions not being change-action differentiable for silly reasons. 

\begin{eg} \label{Eg not all monotone functions are differentiable in PoSet}
Not all $\PoSet$ maps are differentiable in $\PoSet$. To see this explicitly one can insert a constant portion to the translation $T_1$ above as follows:  

\[g : \mN \to \mN ; x \mapsto	\begin{cases}
 							x+1 & x \leq k-1 \\
		 					k 	& k \leq x \leq \ell \\
 							k + x - \ell +1  & \ell +1 \leq x
 							\end{cases}.
 \]
Then by assuming a change-action derivative, $\partial g$, exists and considering the consequences of \textbf{CAD1} and \textbf{CAD2} one sees that such a function is not monotone by evaluating at $(k-1,1)$ and $(k-2,1)$. 
\end{eg}
Although $\partial g$ is not monotone, it is a function which satisfies the commuting diagrams of $\mathbf{CAD1}$ and $\mathbf{CAD2}$ in $\Set$.
Forgetting the order structure on $\PoSet$ also preserves the monoid actions so $g$ is change-action differentiable in $\Set$, but not $\PoSet$.
The main functor we consider later for persistent homology, tracking the lifespans of homology classes, has similar behaviour.
Instead of forgetting structure and working in $\Set$ to resolve the issue, in Section \ref{S CAD in multipersistence} we move to a more general setting: the category of (small) categories and functors, $\Cat$

%SECTION 2

\section{Finite Differences}\label{S2 Finite Differences}

All groups in this section are assumed to be abelian.
In \cite[Section 5.2]{AP&L}, the calculus of finite differences is described in terms of functions $f: G \to H$ between the underlying sets of two abelian groups, $G $ and $H$.
The change actions on the groups are translations given by their respective group operations.
Using the abelian structure on $H$ along with \textbf{CAD1},  one obtains a formula,
\[ \partial f (a,b) = f(a+b) - f(a) ,\]
for the derivative of a given function $f: G \to H$.
This example inspired Alvarez-Picalo and Lemay to introduce Cartesian Difference Categories in \cite{AP&L} as the appropriate context for a flavour of abstract differentiation that captures finite difference operators as derivatives.
In this section, we view the calculus of finite differences in the context of categories and functors by describing a canonical choice of change-action derivatives for functors $F:\mathbb{X} \to \mathbf{Arr}(G)$ where $\mathbb{X}$ has a change-action structure, $G$ is an Abelian group, and $\mathbf{Arr}(G)$ is the arrow category of the `delooping' of $G$, $\mathbf{B}G$. 

Every (Abelian) group $G$ can be viewed as a category, $\mathbf{B}G$, with a single object, $*$, whose morphisms are the group elements and whose composition is given by the group operation.
The arrow category of $\mathbf{B}G$, $\mathbf{Arr}(G)$, has morphisms in $\mathbf{B}G$ as its objects, which are the group elements of $G$, and commuting squares in $\mathbf{B} G$ as its morphisms.
That is, $(f,g) : a \to b$ is a morphism in $\mathbf{Arr}(G)$ whenever the square
\begin{center}
\begin{tikzcd}
	* \dar["a"'] \rar["f"] & * \dar["b"] \\
	* \rar["g"'] & * 
\end{tikzcd}.	
\end{center}
commutes in $\mathbf{B}G$. In other words $fh = ag$. 

When $G$ is abelian, the group operation induces a canonical functor 
\[ + : \mathbf{Arr}(G) \times \mathbf{Arr}(G) \to \mathbf{Arr}(G) \]
defined on objects $ a , b \in G$ by 
\[ + ((a, c)) = a + c . \]

On arrows in $\mathbf{Arr}(G)$ it is defined similarly
\begin{center}
\begin{tikzcd} 
	* \dar["a"'] \rar["f"] & * \dar["b"] \\
	* \rar["g"'] & * 
\end{tikzcd} $\quad  + \quad $
\begin{tikzcd}
	* \dar["c"'] \rar["h"] & * \dar["d"] \\
	* \rar["k"'] & * 
\end{tikzcd}$ \quad = \quad $
\begin{tikzcd}
	* \dar["a+c"'] \rar["f+h"] & * \dar["b+d"] \\
	* \rar["g+k"'] & * 
\end{tikzcd}.
\end{center}

We write $0$ for the functor that picks out the object $0 \in \mathbf{Arr}(G)$ and its identity morphism $(0,0)$.
More precisely if $\mathbf{1}$ is the category with a single object, $*$, and a single morphism, $\id_*$, then the functor we are talking about is 
\[ 0  : \mathbf{1} \to \mathbf{Arr}(G) \ ; \  * \mapsto 0 .\]
The functorial calculus of finite differences involves functors whose codomains are categorified abelian groups, $\mathbf{Arr}(G)$ for some abelian group $G$.
The functors $+$ and $0$ defined above allow the monoid structure on $\mathbf{Arr}(G)$ to be described by the following Lemma.
This monoid acts on itself to provide a change action structure, and we use this in Section \ref{S CAD in multipersistence} to capture the notion of a persistence diagram in a general setting for persistent homology.

\begin{lem}\label{Lem G1 monoid under addition} 
The triple ($\mathbf{Arr}(G)$, $+_\mbG$, $0_\mbG$) is a monoid in $\Cat$  
\end{lem}
\begin{proof}
Let $! : {\mathbf{Arr}(G)}\to \mathbf{1}$ be the functor into the terminal category. Then the following diagrams commute in $\Cat$: 

\begin{enumerate}
\item \[
\begin{tikzcd}[column sep = large, row sep = large]
\mathbf{Arr}(G) \rar["\langle 1_ \mbG {,} ! 0\rangle "] \ar[dr, equals] \dar["\langle ! 0 {,} 1_ \mbG\rangle"'] & \mathbf{Arr}(G)  \times \mathbf{Arr}(G) \dar["+_\mbG"] \\
\mathbf{Arr}(G) \times \mathbf{Arr}(G) \rar["+_\mbG"'] & \mathbf{Arr}(G)
\end{tikzcd}\]
\item \[\begin{tikzcd}[column sep = large, row sep = large]
\mathbf{Arr}(G) \times \mathbf{Arr}(G) \times \mathbf{Arr}(G) \rar["1_ \mbG \times +_\mbG"]\dar["+_\mbG \times 1_ \mbG"'] & \mathbf{Arr}(G) \times \mathbf{Arr}(G) \dar["+_\mbG"]\\
\mathbf{Arr}(G) \times \mathbf{Arr}(G) \rar["+_\mbG"'] & \mathbf{Arr}(G) 	
\end{tikzcd}.
\] 
\end{enumerate} 
The identity triangles on the left commute because $0$ is the identity element in $G$.
The associativity square on the right commutes because the operation in $G$ is associative. 
\end{proof}

The existence of all inverses in $G$ induces another functor,
\[ -_\mbG : \mathbf{Arr}(G) \times \mathbf{Arr}(G) \to \mathbf{Arr}(G) ,\]
defined similarly on a pair of commuting squares by taking differences instead of sums.
More precisely, for any $a, b \in G$, writing $a - b$ in place of $ -(a,b)$ once again, we have 
\begin{center}
\begin{tikzcd} 
	* \dar["a"'] \rar["f"] & * \dar["b"] \\
	* \rar["g"'] & * 
\end{tikzcd} \quad $-$ \quad 
\begin{tikzcd}
	* \dar["c"'] \rar["h"] & * \dar["d"] \\
	* \rar["k"'] & * 
\end{tikzcd} \quad $=$ \quad
\begin{tikzcd} 
	* \dar["a-c"'] \rar["f-h"] & * \dar["b-d"] \\
	* \rar["g-k"'] & * 
\end{tikzcd}.
\end{center}

The difference functor defines a monoid action of $\mathbf{Arr}(G)$ on itself.

\begin{prop}\label{Prop G1 subtraction change action in Cat}
$(\mathbf{Arr}(G), \mathbf{Arr}(G), -_\mbG, +_\mbG, 0_\mbG)$ is a monoid action in $\Cat$. 
\end{prop}
\begin{proof} The diagrams 

\begin{enumerate}
\item \[\begin{tikzcd}[column sep = large, row sep = large]
\mathbf{Arr}(G) \rar["\langle 1  {,} !  0 \rangle "] \ar[dr, equals]  & \mathbf{Arr}(G)  \times \mathbf{Arr}(G) \dar["-"] \\
& \mathbf{Arr}(G)
\end{tikzcd}\]\
\item \[ \begin{tikzcd}[column sep = large, row sep = large]
\mathbf{Arr}(G) \times \mathbf{Arr}(G) \times \mathbf{Arr}(G) \rar["1  \times + "]\dar["- \times 1 "'] & \mathbf{Arr}(G) \times \mathbf{Arr}(G) \dar["-"]\\
\mathbf{Arr}(G) \times \mathbf{Arr}(G) \rar["-"'] & \mathbf{Arr}(G) 	
\end{tikzcd}\]
\end{enumerate}
commute in $\Cat$: For all  $a \in G$, $a-0 = a$ shows $(1)$ commutes. Diagram $(2)$ follows from distributivity: 
\[ (a-c) - b = a - (b + c) \]
\end{proof}

The main result of this section is the following theorem.
It says that very functor from a category with a change action structure into the arrow category of an abelian group is (change-action) differentiable and it does most of the heavy lifting to build the change-action derivatives appearing in persistent homology, as seen in the last section of this paper in Theorem \ref{thm extended homological memory functor is differentiable} and Corollary \ref{cor lifespan functor recovered by derivative}.

\begin{thm}[Functorial Calculus of Finite Difference]\label{thm functors into G1 are differentiable}
Let $(\mbA, \Delta \mbA, \oplus_ \mbA, +_ \mbA, 0_ \mbA)$ be a change action in $\Cat$ and $F: \mbA \to \mathbf{Arr}(G)$ a functor. Then $F$ is differentiable in $\Cat$ with respect to $(\mathbf{Arr}(G), \mathbf{Arr}(G), -, +, 0)$
\end{thm}
\begin{proof}
Define 
\[ \partial F  : \mbA \times \Delta \mbA \to \mathbf{Arr}(G) \]
\noi on an arbitrary morphism
\[ \begin{tikzcd}
 	(A, \alpha) \rar[rr, "(f {,} \varphi)"] && (B, \beta)
\end{tikzcd}\]
in $\mbA \times \Delta \mbA$ by mapping it to the morphism
\[\begin{tikzcd}[column sep = large] 
	F(A) - F(A \oplus_\mbA \alpha) \rar[rr, "F(f) -  F(f \oplus_\mbA \varphi)"] && F(B) - F(B \oplus_\mbA \beta)
\end{tikzcd}\]
in $\mathbf{Arr}(G)$. 
The objects and morphisms in $\mathbf{Arr}(G)$ are the arrows and commuting squares of the delooping $\mathbf{B} G$.
This means the arrows $F(f)$ and $F(f \oplus_\mbA \varphi)$ are pairs
\[ F(f) = \bigr( F(f)_0\ , \  F(f)_1 \bigr) \quad \text{ and } \quad F(f \oplus_\mbA \varphi) = \bigr( F(f \oplus_\mbA \varphi)_0 \ , \ F(f \oplus_\mbA \varphi)_1 \bigr) \] 
where $F(f)_i $ and $F(f \oplus_\mbA \varphi)_i$ are group elements in $G$.
The pair defining $\partial F ((f, \varphi))$ is then 
\[ F(f) - F(f \oplus_\mbA \varphi) = \big( F(f)_0 - F(f \oplus_\mbA \varphi)_0 , F(f)_1  - F(f \oplus_\mbA \varphi)_1 \big)\]
such that the diagram
\begin{center}
\begin{tikzcd}[column sep = huge, row sep = huge]
	* \dar["F(f)_0 - F(f \oplus_\mbA \varphi)_0"']\rar[rr, "F(A) - F(A \oplus_\mbA \delta_A)"] && * \dar["F(f)_1 - F(f \oplus_\mbA \varphi)_1"]\\
	* \rar[rr, "F(B) - F(B \oplus_\mbA \delta_B)"'] && * 
\end{tikzcd} 
\end{center}
commutes in {$\mathbf{B}$}$G$.
This square decomposes in terms of the difference change action structure on $\mathbf{Arr}(G)${:}

\begin{enumerate}
\item[(1) \qquad \qquad ]
\begin{center}
\begin{tikzcd}[column sep = huge, row sep = huge]
	* \dar["F(f)_0"']\rar["F(A)"] & * \dar["F(f)_1"]\\
	* \rar["F(B) "'] & * 
\end{tikzcd}\quad $-$ \quad 
\begin{tikzcd}[column sep = huge, row sep = huge]
	* \dar["F(f \oplus_\mbA \varphi)_0"']\rar[" F(A \oplus_\mbA \delta_A)"] & * \dar["F(f \oplus_\mbA \varphi)_1"]\\
	* \rar["F(B \oplus_\mbA \delta_B)"'] & * 
\end{tikzcd}.
\end{center}
\end{enumerate}

To verify that \textbf{CAD1} and \textbf{CAD2} hold, take any category $\mbX$ and any functors $\chi: \mbX \to \mbA$ and $\delta , \gamma: \mbX \to \Delta \mbA$.
Let $x : X \to Y$ be an arbitrary arrow in $\mbX$ and{,} to simplify the notation in the calculations below{,} relabel $\chi(X) = A$, $\chi(X) = B$, $\delta_A = \delta(X)$, $\delta_B = \delta(Y)$, $\gamma_A =\gamma(X)$, and $\gamma_B = \gamma(Y)$

\[ \begin{tikzcd}[column sep = large]
    A \rar["\chi(x) =: f"] & B
\end{tikzcd} \qquad 
\begin{tikzcd}[column sep = large]
    \delta_A \rar["\delta(x) =: \varphi"] & \delta_B
\end{tikzcd} \qquad 
\begin{tikzcd}[column sep = large]
    \gamma_A \rar["\gamma(x) =: \psi"] & \gamma_B
\end{tikzcd}\] 

To check that $\partial F$ satisfies \textbf{CAD1} we need to see that the diagram, 

\begin{center}
\begin{tikzcd}[column sep = huge, row sep = huge]
	\mbX \rar["\langle \chi {,} \delta \rangle"] \dar[ "\langle \chi{,} \langle \chi{,} \delta \rangle \rangle"'] & \mbA \times \Delta \mbA \dar["\oplus_A"] \\
	\mbA \times ( \mbA \times \Delta \mbA) \dar["F \times \partial F"']& \mbA \dar["F"] \\
	\mathbf{Arr}(G) \times \mathbf{Arr}(G) \rar["-"] & \mathbf{Arr}(G)  
 \end{tikzcd},
\end{center}
commutes in $\Cat$.
The following diagrams in $\mathbf{Arr}(G)$ should be read from top to bottom; the left side of the diagram above is given by
\begin{center}
\begin{tikzcd}[column sep = huge, row sep = huge]
	* \dar["F(f)_0"']\rar["F(A)"] & * \dar["F(f)_1"]\\
	* \rar["F(B) "'] & * 
\end{tikzcd} \quad $-$\quad 
\begin{tikzcd}[column sep = huge, row sep = huge]
	* \dar["\partial F(f{,} \varphi)_0"']\rar["\partial F(A{,} \delta_A)"] & * \dar["\partial F(f{,} \varphi)_1"]\\
	* \rar["\partial F(B{,} \delta_B) "'] & * 
\end{tikzcd}
\end{center}
Expanding the square on the left by (1) above and using the fact that $a-(a-b) = b$ for all $a,b \in G$ gives 
\[
\begin{tikzcd}[column sep = huge, row sep = huge]
	* \dar["F(f \oplus_\mbA \varphi)_0"']\rar[" F(A \oplus_\mbA \delta_A)"] & * \dar["F(f \oplus_\mbA \varphi)_1"]\\
	* \rar["F(B \oplus_\mbA \delta_B)"'] & * 
\end{tikzcd}, \]
which is precisely the other side of the the original diagram in $\Cat$.
This shows the diagram for $\textbf{CAD1}$ commutes. 

To see the first part of \textbf{CAD2}, we need to show that the diagram
\begin{center}
\begin{tikzcd}[column sep = huge, row sep = huge]
 \mbX \rar["\langle \chi{,} !_\mbX 0_{\Delta \mbA} \rangle"] \dar["!_\mbX"'] & \mbA \times \Delta \mbA \dar["\partial F"] \\
 \mathbf{1} \rar["0"'] & \mathbf{Arr}(G)	
 \end{tikzcd},
\end{center}
commutes in $\Cat$.
Expanding the definition of $\partial F$ and noting the identity of the monoid $(\Delta \mbA, +_\mbA, 0_\mbA)$ acts trivially on $\mbA$ quickly shows that 
\begin{align*}
\begin{tikzcd}[ampersand replacement = \&, column sep = large, row sep = large ]
	* \dar["\partial F(f{,} 0)_0"']\rar["\partial F(A{,} 0)"] \& * \dar["\partial F(f{,} 0)_1"]\\
	* \rar["\partial F(B{,} 0) "'] \& * 
\end{tikzcd} 	\qquad 
&= \qquad 
\begin{tikzcd}[ampersand replacement = \&,column sep = large, row sep = large ]
	* \dar["F(f)_0 - F(f)_0"']\rar["F(A)-F(A)"] \& * \dar["F(f)_1-F(f)_1"]\\
	* \rar["F(B) -F(B)"'] \& * 
\end{tikzcd} 	\\
& = \qquad \qquad \qquad \ 
\begin{tikzcd}[ampersand replacement = \&,column sep = large, row sep = large ]
	* \dar["0"']\rar["0"] \& * \dar["0"]\\
	* \rar["0"'] \& * 
\end{tikzcd}	
\end{align*}

This last square is precisely the arrow $(0,0) : 0 \to 0$ in $\mathbf{Arr}(G)$,  which is the image of the composite of the bottom and left functors in the original diagram.
This shows the first part of \textbf{CAD2} holds.
For the second part of \textbf{CAD2} we need to see the diagram,
\begin{center}
\begin{tikzcd}[column sep = huge, row sep = huge]
 \mbX \dar["\langle \langle \chi {,} \delta \rangle {,} \langle \langle \chi{,}\delta \rangle \oplus_A {,} \gamma \rangle \rangle"'] \rar["\langle \chi {,} \langle \delta{,} \gamma \rangle +_\mbA \rangle"] & \mbA \times \Delta \mbA \dar[dd, "\partial F"] \\
 (\mbA \times \Delta \mbA) \times (\mbA \times \Delta \mbA) \dar["\partial F \times \partial F"'] & \\
\mathbf{Arr}(G) \times \mathbf{Arr}(G) \rar["+_\mbG"'] & \mathbf{Arr}(G) 
\end{tikzcd}, 
\end{center}
commutes in $\Cat$.
For this we expand
\begin{center}
\begin{tikzcd}[column sep = huge, row sep = huge]
	* \dar["\partial F(f{,} \varphi +_\mbA \psi)_0"']\rar["\partial F(A{,} \delta_A +_\mbA \gamma_A)"] & * \dar["\partial F(f{,} \varphi +_\mbA \psi)_1"]\\
	* \rar["\partial F(B{,} \delta_A +_\mbA \gamma_A) "'] & * 
\end{tikzcd} 
\end{center}
using the description of $\partial F$ in (1) above, use the fact that the action $\oplus_\mbA$ is compatible with the monoid operation $+_\mbA$, and then add and subtract the following square 
\begin{center}
\begin{tikzcd}[column sep = huge, row sep = huge]
	* \dar["F(f \oplus_\mbA \varphi)_0"']\rar["F(A \oplus_\mbA \delta_A)"] & * \dar["F(f \oplus_\mbA \varphi)_1"]\\
	* \rar["F(B \oplus_\mbA \delta_B) "'] & * 
\end{tikzcd}
\end{center}
using the sum and difference functors on $\mathbf{Arr}(G)$. The result is equal to 
\begin{center}
\begin{tikzcd}[column sep = huge, row sep = huge]
	* \dar["\partial F(f{,} \varphi)_0"']\rar["\partial F(A{,} \delta_A)"] & * \dar["\partial F(f{,} \varphi)_1"]\\
	* \rar["\partial F(B{,} \delta_A ) "'] & * 
\end{tikzcd} \quad $+$ \quad
\begin{tikzcd}[column sep = huge, row sep = huge]
	* \dar["\partial F(f \oplus_\mbA \varphi {,} \psi)_0"']\rar["\partial F(A \oplus_\mbA \delta_A {,}  \gamma_A)"] & * \dar["\partial F(f \oplus_\mbA \varphi {,} \psi)_1"]\\
	* \rar["\partial F(B \oplus_\mbA \delta_A {,} \gamma_A) "'] & * 
\end{tikzcd}. 	
\end{center}
This shows the second part of \textbf{CAD2} holds and proves the theorem. 
\end{proof}

In Section \ref{S CAD in multipersistence} we use Theorem \ref{thm functors into G1 are differentiable} to obtain a change-action derivative that captures a generalized persistence diagram.
To relate this more explicitly to the finite difference calculus in \cite[Section 5.2]{AP&L}, we need the inverse functor, 
\[ \operatorname{inv}: \mathbf{Arr}(G) \to \mathbf{Arr}(G), \]\

\noi defined by sending every element of $G$ in sight to its additive inverse: 

\[ \operatorname{inv} \left( \begin{tikzcd}
	* \dar["f"']\rar["a"] & * \dar["g"]\\
	* \rar["b "'] & * 
\end{tikzcd} \right) \quad = \quad 
\begin{tikzcd}
	* \dar["-f"']\rar["-a"] & * \dar["-g"]\\
	* \rar["-b"'] & * 
\end{tikzcd}\] 

\noi More precisely, there is another change-action $(\mathbf{Arr}(G), \mathbf{Arr}(G), +, +, 0)$ which is used to describe the calculus of finite differences in \cite[Section 5.2]{AP&L}.
Similarly to how $\partial F$ was defined in Theorem \ref{thm functors into G1 are differentiable}, any functor $F: \mbA \to \mathbf{Arr}(G)$ is differentiable with respect to this change action structure, $(\mathbf{Arr}(G), \mathbf{Arr}(G), +, +, 0)$.
Its derivative, denoted $- \partial F$ here, is given by post composing $\partial F$ with the inversion functor.

\begin{prop}\label{Prop functors into G1 are differentiable in a different way}
Any functor $F : \mbA \to \mathbf{Arr}(G)$ is differentiable with respect to the change action $(\mathbf{Arr}(G), \mathbf{Arr}(G), +, +, 0)$ and its derivative is given by the composite:
\begin{center}
\begin{tikzcd}[column sep = large, row sep = large]
	\mbA \times \Delta \mbA \rar["\partial F"] \ar[dr, "-\partial F"']  & \mathbf{Arr}(G) \dar["inv"] \\
	 & \mathbf{Arr}(G)
\end{tikzcd}
\end{center}
\end{prop}
\begin{proof}
By definition, 
\[ -\partial F  : \mbA \times \Delta \mbA \to \mathbf{Arr}(G) \]
sends an arrow 
\[ \begin{tikzcd}[column sep = huge]
 	(A, \delta_A) \rar["(f {,} \varphi)"] &(B, \delta_B)
\end{tikzcd}\]
in $\mbA \times \Delta A$ to the arrow, 
\[
% \begin{tikzcd}
%	\ \rar[mapsto , "-\partial F"] & \
%\end{tikzcd} 
\begin{tikzcd}
	F(A \oplus_\mbA \delta_A) - F(A) \rar[rr, " F(f \oplus_\mbA \varphi) - F(f)"] 
	 &&  F(B \oplus_\mbA \delta_B) - F(B) 
\end{tikzcd}
\]
in $\mathbf{Arr}(G)$, which is precisely how the composite $\operatorname{inv} \circ \partial F$ is defined.
The axioms \textbf{CAD1} and \textbf{CAD2} follow from Theorem \ref{thm functors into G1 are differentiable} and changing signs. 
\end{proof}

%SECTION 3

\section{Multipersistence}\label{S3 Multipersistence}

For a sufficiently finite one-parameter filtration, the persistent homology can be extracted from the rank of the pair group using the inclusion-exclusion principle \cite[Section 2, Tame Functions]{E&H}.
In this section we consider a generalization of the pair group for sufficiently finite multi-parameter filtrations which was originally given by McCleary and Patel in a paper which has since been retracted.
The main area of contention in their work seemed to be the generalized inclusion-exclusion principle used to count the rank of the generalized pair-group, as opposed to the pair group itself.
The purpose of this paper is not to address this counting problem, but to show that the rank of the pair group can be described in terms of a change-action derivative. 

We work directly with filtrations of chain complexes indexed by a partially ordered set.
The partially ordered set is viewed as a topological space using the upward closed subsets so we can consider cycles and boundaries that appear at or before the boundary of a given open.
For a pair of nested opens $U \subseteq V$, the intersection of `cycles that appear by $U$' and `boundaries that appear by $V$' represents the homology cycles who lived during or before the region $V - U$.
Identifying those cycles that were born strictly before $U$ or that died strictly before $V$ leaves us with the cycles whose lifespan is precisely the region $V - U$. 

\subsection{Filtrations}\label{SS Filtrations}\

Let $\cA$ be an abelian category, $\cP$ a partially ordered set, and $\mathbf{Ch}(\cA)$ the category of chain complexes in $\cA$.
Since homology cycles with short lifespans are generally considered `noise,' and arbitrary filtrations can get rather `wild' in general, we are forced to restrict our attention to sufficiently `tame' filtrations. 

For example when $\cP = \mR$ one can construct a filtration whose persistent homology changes at every irrational number and the persistence diagram is therefore riddled with `noise.' 

%Do I want to identify subobjects or just take the subcategory of monos in the slice? 
%and $\mathbf{Sub}_{\cA}(X)$ the category of subobjects for any object $X$ in $\cA$. 

\begin{defn} \label{Def filtration}
Let, \begin{tikzcd}	\cP \rar[r, "F"] & \mathbf{Ch}(\cA)\end{tikzcd} be a functor. $F$ is

\begin{center}
\begin{itemize}
\item a \textit{filtration} if for any $x \leq y \in \cP$, $F(x \leq y) = F(X) \rightarrowtail F(Y)$ is monic in $\ChA$. \\
\item \textit{bounded} if the colimit, denoted $F^\fty$, and limit, denoted $F^{-\fty}$, both exist.
\item \textit{finitely supported} if and only if there exists a finite sub-poset $\cP' \subseteq \cP$ such that $F(x) \cong F(y)$ whenever there is no $z \in \cP'$ such that $x \leq z < y$. 
\item a \textit{tame filtration} if it is a finitely supported and bounded filtration
\end{itemize}.
\end{center}
\end{defn}

Working with tame filtrations allows one to restrict $F$ to a finite poset $\cP'$ without losing information about the homology and realistically count the homology cycles throughout the poset.
For applications of persistent homology to data analysis, one usually deals with tame filtrations where $\cP = \mR^n$ and $\cA = \mathbf{Vec}^{fd}_K$ is the category of finite-dimensional vector spaces over a fixed field, $K$, usually with charactestic zero or $K = \mZ/2\mZ$. 

For each $x \in \cP$ we write $F(x)$ for the chain complex 
\[\begin{tikzcd}
... \rar["\partial_{n+2}^x"] & F_{n+1}(x) \rar["\partial_{n+1}^x"] & F_n(x) \rar["\partial_n^x"] & F_{n-1}(x) \rar["\partial_{n-1}^x"] & ...
\end{tikzcd}.
\]

\subsection{Subobjects of Cycles}
\label{SS Subobjects of Cycles}

There are two important subobjects to consider in each degree for chain complexes when discussing homology. 

\begin{defn}\label{Def birth and death subobjects}
For each $n \in \mathbb{Z}$, $x \in \cP$, let
\[ \cZ F_n (x) := \ker \partial_n^x \quad \text{ , } \quad  \begin{tikzcd} \cZ F_n(x) \rar[tail, "\iota_n^x"]  &F_n(x)\end{tikzcd} \]
and 
 \[ \cB F_n(x) := \text{Im}\partial_{n+1}^x \quad \text{ , }  \begin{tikzcd} \cB F_n(x) \rar[tail, "\text{im}\partial_{n+1}^x"] & F_n(x)\end{tikzcd} \]
\end{defn}

\begin{rems}\label{subobjects of cycle colimit}
Both $\iota_n^x$ and im$\partial_{n+1}^x$ are subobjects of $\cZ F_n^\fty$.
Their domains, $\cZ F_n(x)$ and $\cB F_n(x)$, represent the `$n$-cycles,' and `boundary $n$-cycles,' of $F$ at $x$ respectively.
Persistent homology wants to know at which points in $\cP$ a class of $n$-cycles appears and at which (possibly later) points are they boundaries of a higher cell. 
\end{rems}

\begin{defn}\label{defn ZBFn on points} 
Let $\Sub_\cA( \cZ F_n^\fty)$ denote the category whose objects are equivalence classes of monomorphisms
\[ \begin{tikzcd}
 A \rar[tail, "\alpha"] & \cZ F_n^\fty 
 \end{tikzcd}
\]
\noi in $\cA$ and whose morphisms, $f : \alpha \to \beta$, are given by commuting triangles, 
\[ \begin{tikzcd}
 A \dar[rr, "f"] \ar[dr, tail, "\alpha"'] && B \ar[dl, tail, "\beta"]\\
 & \cZ F_n^\fty	&  
 \end{tikzcd},
\]
\noi in $\cA$. 
\end{defn}\

\noi Since $\cA$ has finite limits, for any two pairs of points $x, y \in \cP$ we can take the following pullback in $\cA$. 
\begin{center}
\begin{tikzcd}[column sep = huge, row sep = huge]
\cZ F_n(x) \cap \cB F_n(y) \arrow[dr, phantom, "\usebox\pullback" , very near start, color=black] \rar[tail] \dar[tail] & \cB F_n(y) \dar[tail] \\
\cZ F_n(x) \rar[tail] &\cZ F_n^\fty
\end{tikzcd}
\end{center}
This pullback is a product in $\Sub_\cA(\cZ F_n^\fty)$ and represents the cycles of $F$ that appear by $x$ and become boundaries by $y$. 

\subsection{Open Lattices of Up-Sets and Their Arrow Categories}\label{SS Upper Sets}

The intersection, $\cZ F_n(x) \cap \cB F_n(y)$, gives us some data about the persistent homology of $F$, but to compare each of these intersections for all pairs of elements $x \leq y$ in $\cP$ can be computationally infeasible.
In this section we extend Definitions \ref{Def birth and death subobjects} and \ref{defn ZBFn on points} to pairs of nested open subsets in a suitable topology on $\cP$.

\begin{defn}\label{Def upper-sets}
For any poset $\cP$, let $\tau$ be the set of all $U \subseteq \cP$ such that if $x \in U \text{ and } x \leq y$ then $y \in U$.
The elements of $\tau$ are called the \textit{up-sets} of $\cP$.
\end{defn}

\begin{prop} \label{upper sets form a topology}
The poset $\tau$ is a topology which is closed under arbitrary intersections. 
\end{prop}
\begin{proof} Immediate from definition. 
%To see $\tau$ is a topology we need to show that $\cP$ and $\emptyset$ are in $\tau$, that $\tau$ is closed under arbitrary unions, and that $\tau$ is closed under finite intersections. We'll show it's actually closed under arbitrary intersections as well.  

%Suppose $x \in \cP$ and $ x \leq y$. Well $y$ must be in $\cP$ or else we couldn't have assumed $x \leq y$ since the relation $\leq$ is defined on $\cP$. Hence $\cP \in \tau.$ Similarly, $x \in \emptyset$ and $x \leq y$ implies $ y \in \emptyset$ for vacuous reasons because there's no such $x$ in the first place.

%Now suppose $G_k \in \tau$ for all $k \in K$ some arbitray index set. If $x \in \cup_{k \in K} G_k $ and $x \leq y$, then $x \in G_k$ for some $k$ and $x \leq y$ implies that $y \in G_k$ for that same $k$. Hence $y \in \cup_{k \in K} G_k$, and $\cup_{k \in K} G_k \in \tau$. 

%To see closure under arbitrary intersections, let $G_k \in \tau$ for all $k \in K$ some arbitrary indexing set, $x \in \cap_{k \in K} G_k$ and $x \leq y$. Then $x \in G_k$ and $x \leq y$ for each $k \in K$ implies that $y \in G_k$ for each $k \in K$ and therefore $y \in \cap_{k \in K} G_k$. This shows $\cup_{k \in K} G_k \in \tau$.  
\end{proof}

The up-set topology, $\tau$, is a poset itself under set inclusion and can be viewed as a category whose objects are the open subsets and whose morphisms are the inclusions.
The opposite category, $\tau^{op}$, is more relevant for our purposes.
The morphisms of $\tau^{op}$ are called restrictions.
To be precise about where exactly cycles are born and where exactly they become boundaries, we need to talk about maximal proper subobjects in $\tau^{op}$.
These are minimal singleton covers of opens in $\tau$ so we call them \emph{blankets}.

%We use the term `cover' to mean a maximal proper subobject as Johnstone does in \cite[page 19, A.1.3]{sketches}.
%\biggeoff{I thought cover to Johnstone meant a map $f:U \to V$ in a poset $X$ for which $f$ could \emph{not} factor through any proper subobject of $V$ and not that covers were maximal proper subobjects? Or I'm getting cofused about the $\tau$ versus $\tau^{op}$ thing... Maybe saying ``$\tau$-cover'' or ``$\tau^{op}$-cover'' would help make this distinction clear?}
%\bigdeni{Yeah I agree this is confusing I don't like this word for this at all, a maximal proper subobject in $\tau^{op}$ is a minimal singleton cover in $\tau$, I call it a blanket now  }

\begin{defn}\label{Def blanket}
	A restriction $\begin{tikzcd} U \rar["\supset"] & V \end{tikzcd}$ in $\tau^{op}$ is a \textit{blanket} if %it is a maximal proper subobject. That is, 

	\[\begin{tikzcd}W \rar["\supseteq"] & U \rar["\supset"] &  V \end{tikzcd} \quad \text{ implies } \quad \begin{tikzcd} W \rar[equals] & U \end{tikzcd}. \]\
\end{defn}\

\begin{eg}
Consider $\mR$ as a poset with its usual ordering.
For each $x \in \mR$ the upper set principally generated by $x$ is the left-closed interval 
\[ \{ y \in \mR : x \leq y \} = [x, \fty).  \]
Any two upper sets are comparable since $\mR$ is totally ordered.
Every upward closed subset is an intersection of principally generated upward closed sets, for example 
\[ (a, \infty) = \bigcap_{n =1}^\fty [a + 1/n, \fty)\]
Since the infimums of every upward closed subset are comparable in $\mR$, every pair of upward closed subsets are comparable in $\tau$, and hence every pair of upward closed subsets are comparable in $\tau^{op}$.
The only blankets in $\tau^{op}$ are the restrictions of a left-closed interval to the open interval that deletes its left endpoint:  
\[ \begin{tikzcd} \lbrack x, \fty) \rar["\supset"] &  (x, \fty)\end{tikzcd} \]
\end{eg} 

\begin{eg}\label{Egs upper-sets in R^2}
The sub-poset $ \cP := [0,\fty) \times [0, \fty) \subset \mR^2$ with the product partial order, where $(a,b) \leq (a',b')$ if and only if $a \leq a'$ \emph{and} $b \leq b'$, is convenient for visualizing the complexity with multiple parameters.
The up-set generated by an element $z \in \cP$ is defined to be $[z, \infty) = \{ w \in \cP : z \leq w \}$.
The interval notation here is purely formal and the up-set is really a rectangular region in the plane as pictured in pink below.  

\begin{figure}[!h]
\begin{tikzpicture}
\draw [step=1cm,white,very thin] (-1.9,-1.9) grid (4,4);
\draw[thick,->] (-1,-1) -- (3.5,-1) node[anchor=north west] {};
\draw[thick,->] (-1,-1) -- (-1,3.5) node[anchor=south east] {};

\begin{scope}[transparency group]
\begin{scope}[blend mode = normal, fill opacity = 20]
%\fill[magenta!5] (1,0) rectangle (6,6);
%\fill[magenta!5] (0,1) rectangle (6,6);
%\fill[magenta!10] (1,1) rectangle (6,6);
\fill[magenta!15] (0,0) rectangle (4,4);
%\fill[magenta!20] (4,4) rectangle (6,6);
\end{scope}
\end{scope}

\draw[thick, dotted] (0,0) -- (0,4);
\draw[thick, dotted] (0,0) -- (4,0);

%nodes
\draw (0,0) circle(0.1 cm) node[anchor=north east] {$z$} ;
\end{tikzpicture}
\caption{The up-set, $[z, \infty)$, generated by $z \in \mR^2$}
\end{figure}
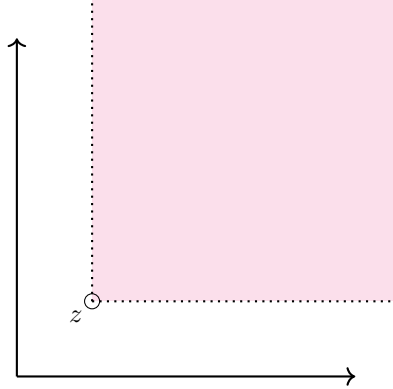 
An example of a subset of mutually incomparable elements in $\cP = \mR^2$ is the set of elements $\{ x,y \}$ in Figure \ref{Fig2 incomparable upsets and their colimit}. 
One way to see they're incomparable is that the unique line connecting $x$ to $y$ in $\mR^2$ has a negative slope.
The set-theoretic union of the up-sets they generate is their limit/meet in $\tau^{op}$, pictured in pink, and their set-theoretic intersection is their colimit/join in $\tau^{op}$, pictured in teal area.
We'll discuss blankets shortly after in Example \ref{Eg blankets of up-sets in R^2}.

\begin{figure}[h!]
\begin{tikzpicture}
\draw [step=1cm,white,very thin] (-1.9,-1.9) grid (4,4);
\begin{scope}[transparency group]
\begin{scope}[blend mode = normal, fill opacity = 20]
%\fill[magenta!5] (1,0) rectangle (6,6);
%\fill[magenta!5] (0,1) rectangle (6,6);
\fill[magenta!10] (0,1) rectangle (4,4);
\fill[magenta!10, overlay] (1,0) rectangle (4,4);
%\fill[magenta!25, overlay] (1,1) rectangle (4,4);
%\fill[magenta!10] (-1,1) --  (-1,4)  -- (4,4) --  (4,-1) -- (1,-1) -- cycle;
%\fill[magenta!20] (4,4) rectangle (6,6);
%\clip[insert path = {(-1,-1) rectangle (4,4)} ] (-1,-1) -- ++ (-1,2) -- ++ (2,-1) -- (-1,-1) cycle ;
\end{scope}
\end{scope}
\draw[thick,->] (-1,-1) -- (3.5,-1) node[anchor=north west] {};
\draw[thick,->] (-1,-1) -- (-1,3.5) node[anchor=south east] {};
%ray
\draw[thick ,dotted] (0,1) -- (0,4) node[anchor=south east] {};
\draw[thick ,dotted] (0,1) -- (4,1) node[anchor=south east] {};
\draw[thick ,dotted] (1,0) -- (4,0) node[anchor=south east] {};
\draw[thick ,dotted] (1,0) -- (1,4) node[anchor=south east] {};
%\draw[thick,red] (-1,1) -- (1,-1) node[anchor=south east] {};
%nodes
\draw (0,1) circle(0.1 cm) node[anchor=north east] {$x$} ;
\draw (1,0) circle(0.1 cm) node[anchor=north east] {$y$} ;
\end{tikzpicture} 
\qquad \qquad 
\begin{tikzpicture}
\draw [step=1cm,white,very thin] (-1.9,-1.9) grid (4,4);
\begin{scope}[transparency group]
\begin{scope}[blend mode = normal, fill opacity = 20]
%\fill[magenta!5] (1,0) rectangle (6,6);
%\fill[magenta!5] (0,1) rectangle (6,6);
%\fill[magenta!10] (0,1) rectangle (4,4);
%\fill[magenta!10, overlay] (1,0) rectangle (4,4);
\fill[teal!10, overlay] (1,1) rectangle (4,4);
%\fill[magenta!10] (-1,1) --  (-1,4)  -- (4,4) --  (4,-1) -- (1,-1) -- cycle;
%\fill[magenta!20] (4,4) rectangle (6,6);
%\clip[insert path = {(-1,-1) rectangle (4,4)} ] (-1,-1) -- ++ (-1,2) -- ++ (2,-1) -- (-1,-1) cycle ;
\end{scope}
\end{scope}\draw[thick,->] (-1,-1) -- (3.5,-1) node[anchor=north west] {};
\draw[thick,->] (-1,-1) -- (-1,3.5) node[anchor=south east] {};
%ray
\draw[thick ,dotted] (0,1) -- (0,4) node[anchor=south east] {};
\draw[thick ,dotted] (0,1) -- (4,1) node[anchor=south east] {};
\draw[thick ,dotted] (1,0) -- (4,0) node[anchor=south east] {};
\draw[thick ,dotted] (1,0) -- (1,4) node[anchor=south east] {};
%nodes
\draw (0,1) circle(0.1 cm) node[anchor=north east] {$x$} ;
\draw (1,0) circle(0.1 cm) node[anchor=north east] {$y$} ;
\end{tikzpicture}

\caption{Two incomparable up-sets, $[x,\infty)$ and $[y, \infty)$, in $\tau$, their limit/meet in $\tau^{op}$ (pink), and their colimit/join in $\tau^{op}$ (teal).}\label{Fig2 incomparable upsets and their colimit}
\end{figure}
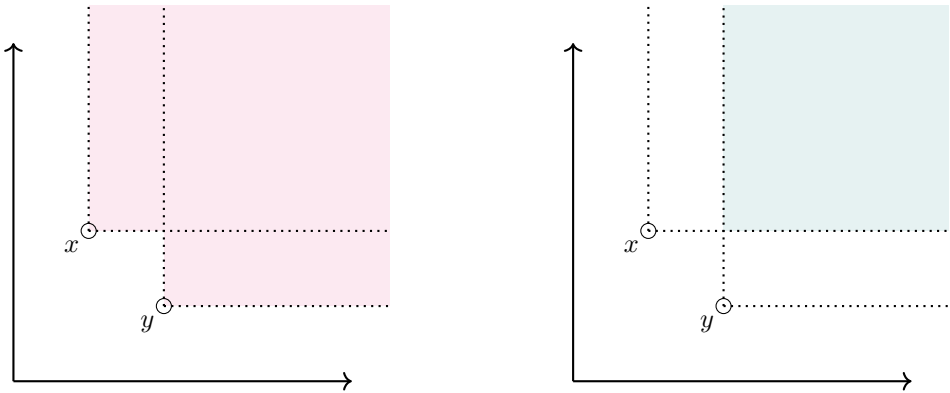
\end{eg}

Let $\mathbf{Arr}(\tau^{op})$ denote the category whose objects are arrows in $\tau^{op}$ and whose morphisms are commuting squares in $\tau^{op}$.
Note that since $\tau^{op}$ is a preoder, any diagram in $\tau^{op}$ is automatically commutative.
In particular, an object in the arrow category is a pair of comparable subsets of $\tau^{op}$ and a morphism in the arrow category witnesses that the domains and codomains are comparable respectively.
The set difference, $U \setminus V$, makes sense for a pair of comparable subsets, $U \supset V$, and a homology cycle of a filtration that is born at the initial boundary of $U$ and dies at the initial boundary of $V$ has a `lifespan' of $U \setminus V$.

Since $\mathbf{Arr}(\tau^{op})$ is a poset category, all of its arrows are monic and a \textit{proper subobject} of $(V_0, V_1)$ is any map \begin{tikzcd} (U_0, U_1) \rar[] & (V_0,V_1) \end{tikzcd} in the arrow category which is not an isomorphism.
This means one or both of $U_0 \supseteq V_0$ or $U_1 \supseteq V_1$ are strict inclusions of open sets.  

\begin{rem}\label{rem characterizing blankets in Arr(tauo^p)}
By Definition \ref{Def blanket}, a morphism \begin{tikzcd}  (U_0, U_1) \rar[-{Triangle[open]}]  & (V_0, V_1), \end{tikzcd} is a blanket in $\mathbf{Arr}(\tau^{op})$ if it is isomorphic to any other proper subobject of $(V_0, V_1)$ it factors through.
So if 

\begin{center}
\begin{tikzcd}
	(U_0, U_1) \ar[dr] \rar[rr, tail, -{Triangle[open]}] && (V_0, V_1) \\
	& (K_0, K_1) \ar[ur, tail] & 
\end{tikzcd},
\end{center}
commutes in $\mathbf{Arr}(\tau^{op})$, and the horizontal map on top is a blanket, the first arrow in the factorization must be an isomorphism, $(U_0, U_1) \cong (K_0, K_1)$.
In particular $U_0 = K_0$ and $U_1 = K_1$ as sets.
This means such a blanket in $\mathbf{Arr}(\tau^{op})$ is a pair of maps ,\begin{tikzcd} U_0 \rar[]  & V_0 \ \end{tikzcd}, \begin{tikzcd} U_1 \rar[]  & V_1 \ \end{tikzcd}, in $\tau^{op}$ such that exactly one of them is a blanket in $\tau^{op}$ and the other is an isomorphism.
\end{rem} 

%\biggeoff{I changed the \texttt{\textbackslash begin\{figure\}} code as before (in this case I just added the `\texttt{!}'. I also moved and uncommented your to the place where it will move the numbers properly and made the forward reference (namely Figure \ref{Fig3 blankets in tauop}) wanted to use work.}
%\bigdeni{thank!}

\begin{eg}\label{Eg blankets of up-sets in R^2}

Consider the first quadrant of $\mR^2$ and let $\cP := \{r_0,r_1,k_0,k_1,k_2 \} \subseteq \mR^2$ be a finite subposet.
Also let $R_i$ and $K_j$ denote the upper sets generated by $r_i$ and $k_j$ for $0 \leq i \leq 1$ and $0 \leq j \leq 2$ respectively, as shown in Figure \ref{Fig3 blankets in tauop} below with varying shades of magenta. 

\begin{figure}[h!]

\begin{tikzpicture}
%\draw (0,0) rectangle (4,4);
%\draw (0,0) parabola (4,4);
%\draw (0,0) .. controls (0,4) and (4,0) .. (4,4);
%\draw[red,thick,dashed] (2,2) circle (3cm);
%\draw (2,2) ellipse (3cm and 1cm);
%\draw (3,0) arc (0:75:3cm);
\draw [step=1cm,white,very thin] (-1.9,-1.9) grid (5.9,5.9);
%\shadedraw[inner color=blue,outer color=red, draw=black] (0,0) rectangle (4,4);
\draw[thick,->] (-1,-1) -- (5.5,-1) node[anchor=north west] {};
\draw[thick,->] (-1,-1) -- (-1,5.5) node[anchor=south east] {};
%\foreach \x in {0,1,2,3,4}
%    \draw (\x cm,1pt) -- (\x cm,-1pt) node[anchor=north] {$\x$};
%\foreach \y in {0,1,2,3,4}
%    \draw (1pt,\y cm) -- (-1pt,\y cm) node[anchor=east] {$\y$};
\begin{scope}[transparency group]
\begin{scope}[blend mode = normal, fill opacity = 20]
\fill[magenta!5] (1,0) rectangle (6,6);
\fill[magenta!5] (0,1) rectangle (6,6);
\fill[magenta!10] (1,1) rectangle (6,6);
\fill[magenta!15] (3,3) rectangle (6,6);
\fill[magenta!20] (4,4) rectangle (6,6);
\end{scope}
\end{scope}

%nodes
\draw (1,1) circle(0.1 cm) node[anchor=south west] {$r_0$} ;
\draw (0,1) circle(0.1 cm) node[anchor=south west] {$k_0$} ;
\draw (1,0) circle(0.1 cm) node[anchor=south west] {$k_2$} ;
\draw (3,3) circle(0.1 cm) node[anchor=south west] {$k_1$} ;
\draw (4,4) circle(0.1 cm) node[anchor=south west] {$r_1$} ;

\draw[thick, dotted, -] (0,1) -- (6,1);
\draw[thick, dotted, -] (0,1) -- (0,6);
\draw[thick, dotted, -] (1,0) -- (1,6);
\draw[thick, dotted, -] (1,0) -- (6,0);
\draw[thick, dotted, -] (3,3) -- (6,3);
\draw[thick, dotted, -] (3,3) -- (3,6);
\draw[thick, dotted, -] (4,4) -- (4,6);
\draw[thick, dotted, -] (4,4) -- (6,4);

\end{tikzpicture}
\caption{An  example of blankets in the arrow category of $\tau^{op}$.}\label{Fig3 blankets in tauop}
\end{figure}
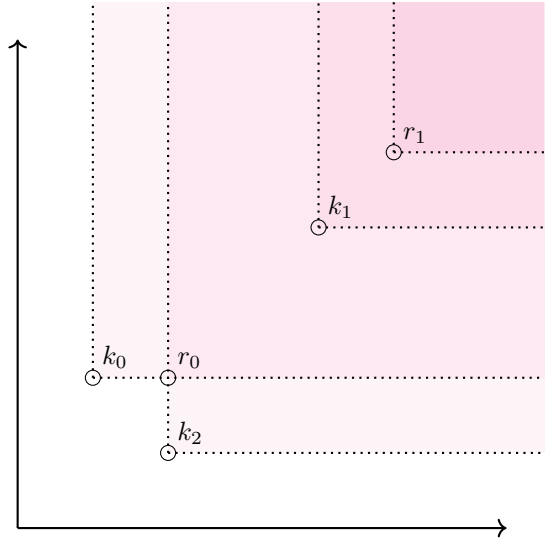

%Notice $K_0,K_2 \supset R_0 \supset K_1 \supset  R_1$ where the left and right restrictions are blankets in $\tau |_{\cP}^{op}$, the topology of upper sets on $\cP$. From now on we'll omit the `restricted to $\cP$' notation when discussing the topologies of relevant sub-posets unless it becomes unclear. 

In this case the blankets of $(R_0, R_1)$ can be obtained by fixing exactly one of $R_0$ or $R_1$ and choosing blankets in $\tau^{op}$ in place of the other.
This implies $(K_0,R_1), (R_0,K_1),$ and $(K_2,R_1)$ are the blankets of $(R_0,R_1)$. 
%\end{tikzcd} 
\end{eg}

\begin{eg}\label{Eg blankets of upper-sets in R^2 may be weird}
One might hope that to find all the blankets of a given object $(U,V)$ in $\mathbf{Arr}(\tau^{op})$ they could fix one of the opens and choosing blankets in $\tau^{op}$ for the other.
This doesn't always work because it may be that some blankets of one open in $\tau{^op}$ are not comparable with the other fixed open.
For example, let $X_i, Y_j$ be the upper-sets generated by $x_i,y_j$ for $i=0,1$ and $j=0,1,2$ respectively in Figure \ref{Fig4 blankets and incomparable} below:

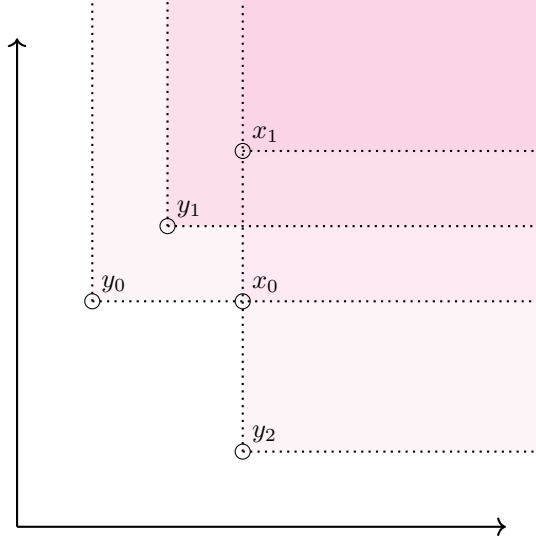
\begin{figure}[h!]
\begin{tikzpicture}
%\draw (0,0) rectangle (4,4);
%\draw (0,0) parabola (4,4);
%\draw (0,0) .. controls (0,4) and (4,0) .. (4,4);
%\draw[red,thick,dashed] (2,2) circle (3cm);
%\draw (2,2) ellipse (3cm and 1cm);
%\draw (3,0) arc (0:75:3cm);
\draw [step=1cm,white,very thin] (-1.9,-1.9) grid (5.9,5.9);
%\shadedraw[inner color=blue,outer color=red, draw=black] (0,0) rectangle (4,4);
\draw[thick,->] (-1,-1) -- (5.5,-1) node[anchor=north west] {};
\draw[thick,->] (-1,-1) -- (-1,5.5) node[anchor=south east] {};
%\foreach \x in {0,1,2,3,4}
%    \draw (\x cm,1pt) -- (\x cm,-1pt) node[anchor=north] {$\x$};
%\foreach \y in {0,1,2,3,4}
%    \draw (1pt,\y cm) -- (-1pt,\y cm) node[anchor=east] {$\y$};
\begin{scope}[transparency group]
\begin{scope}[blend mode = normal, fill opacity = 20]
\fill[magenta!5] (2,0) rectangle (6,6);
\fill[magenta!5] (0,2) rectangle (6,6);
\fill[magenta!10] (2,2) rectangle (6,6);
\fill[magenta!15] (1,3) rectangle (6,6);
\fill[magenta!20] (2,4) rectangle (6,6);
\end{scope}
\end{scope}

%nodes
\draw (2,2) circle(0.1 cm) node[anchor=south west] {$x_0$};
\draw (0,2) circle(0.1 cm) node[anchor=south west] {$y_0$};
\draw (2,0) circle(0.1 cm) node[anchor=south west] {$y_2$};
\draw (1,3) circle(0.1 cm) node[anchor=south west] {$y_1$};
\draw (2,4) circle(0.1 cm) node[anchor=south west] {$x_1$};

\draw[thick, dotted, -] (0,2) -- (6,2);
\draw[thick, dotted, -] (0,2) -- (0,6);
\draw[thick, dotted, -] (2,0) -- (2,6);
\draw[thick, dotted, -] (2,0) -- (6,0);
\draw[thick, dotted, -] (1,3) -- (6,3);
\draw[thick, dotted, -] (1,3) -- (1,6);
%\draw[thick, dotted, -] (2,4) -- (2,6);
\draw[thick, dotted, -] (2,4) -- (6,4);

\end{tikzpicture}

\caption{Comparable opens may have individual blankets in $\tau^{op}$ that are incomparable
%{Geoff suggests a caption for \LaTeX\, numbering fun times high five good job Geoff!}
}\label{Fig4 blankets and incomparable}
\end{figure}\

The blankets of $X_1$ in $\tau^{op}$ are $Y_1$ and $X_0$ and the blankets of $X_0$ are $Y_0$ and $Y_2$. 
In this case, the blankets of $(X_0,X_1)$ in $\mathbf{Arr}(\tau^{op})$ are only $(Y_0, X_1)$ and $ (Y_2,X_1)$. 
This is because $(X_0, Y_1)$ is not an arrow in the arrow category due to the fact that $x_0$ and $y_1$ are incomparable in (the product partial order on) $\mR^2$. 
This is good for persistent homology, because a cycle that appears at $x_0$ cannot become a boundary at $y_1$, it would have to appear at or before $y_1$ in the partial order. 
\end{eg}

%\biggeoff{It might be nice to have an explanation of the purpose of the example below.}

The following example shows how to recover `lifespans' from objects in $\mathbf{Arr}(\tau^{op})$, namely as set-differences of comparable opens in $\tau$.

\begin{eg}
Let $\cP= \{ 0,1,2 \} \subset \mR$ and let $\tau$ be the topology of upper sets on $\cP$. The blanket of $[i, \fty )$ in $\tau^{op}$ is $[i-1, \fty )$ for $i = 1,2$ and $[0,\fty)$ is blanketed by the empty set.
Notice that a pair of comparable upper sets, 
\[ [i, \fty ) \supseteq [j, \fty)\]
determines a `lifespan' 
\[ [i, \fty ) \setminus [j, \fty) = [i , j) \]
which is empty exactly when the upper-sets are equal.
If an $n$-cycle in a one-parameter filtration appears  at time $i$ and becomes a boundary of an $(n+1)$-cycle a at time $j$, then the interval $[i,j)$ is precisely the lifespan of the homology class represented by that $n$-cycle. 
\end{eg}

\subsection{Cycle-Subobjects as Functors }\label{SS Cycle-Subobjects as Functors}

Let $\tau$ be the topology of upper sets on an arbitrary poset as in Definition \ref{Def upper-sets} and let $F$ be a tame filtration which we restrict to a finite sub-poset $\cP$ where the changes in homology take place.
In this section we extend the definitions of $\cZ F_n$ and $\cB F_n$ in Subsection \ref{SS Subobjects of Cycles} to arbitrary opens in order to get two functors $\tau^{op} \to \Sub_\cA(\cZ F_n^\fty)$. 

For any open set $U \subseteq \cP$, view $U$ as a discrete category, i.e.\@ with objects the elements of $U$ and only identity arrows, and consider the diagram:
\[ \begin{tikzcd} U \rar["\cZ F_n "] & \Sub_\cA(\cZ F_n^\fty) \end{tikzcd} \]
defined by sending an object $u$ in $U$ to the subobject
\[ \begin{tikzcd}
    \cZ F_n(u) \rar[tail]  & \cZ F_n^\fty
\end{tikzcd}.\]
Let $\cZ F_n(U)$ denote the limit of this discrete diagram. 
Similarly, for any open $V \subseteq \cP$, let $ \cB F_n(V) $ denote the limit of the diagram 
\[ \begin{tikzcd}
    V \rar["\cB F_n"] & \Sub_\cA(\cZ F_n^\fty)
\end{tikzcd}\]
defined by sending an object $v$ in $V$ to the subobject,  
\[ \begin{tikzcd}
\cB F_n(v) \rar[tail] &  \cB F_n^\fty
\end{tikzcd}. \]
When $F$ is a finitely constructible filtration, these limits are finite and exist because abelian categories have finite limits and because their subobject lattices are finitely complete and cocomplete.
For any opens $U \supseteq V$,  for each $v \in V$ there's a $u \in U$ with $u \leq v$ in $\cP$. 
This means the diagrams

\begin{center}
\begin{tikzcd} 
\cZ F_n(U) \ar[dr, tail] \rar[rr, tail, bend left] \rar[tail]& \cZ F_n(u) \ar[d, tail]\rar[tail] &\cZ F_n(v) \ar[dl, tail]\\
& \cZ F_n^\fty & 
\end{tikzcd} \qquad 
\begin{tikzcd} 
\cB F_n(U) \ar[dr, tail] \rar[rr, tail, bend left] \rar[tail]& \cB F_n(u) \ar[d, tail] \rar[tail] &\cB F_n(v) \ar[dl, tail]\\
& \cZ F_n^\fty &
\end{tikzcd},
\end{center}
commute by definition of limits.
Taking the collection of all of the composites on top, for each $v \in V$, induces the unique horizontal maps
\begin{center} 
\begin{tikzcd}
\cZ F_n(U) \ar[dr, tail] \rar[rr, dashed, tail] && \cZ F_n(V) \ar[dl, tail]  \\
&\cZ F_n^\fty &
\end{tikzcd} 
\qquad 
\begin{tikzcd}  \cB F_n(U) \ar[dr, tail] \rar[rr, dashed, tail] && \cB F_n(V) \ar[dl, tail] \\
&\cZ F_n^\fty &
\end{tikzcd}
\end{center}
by the universal properties of each limit, respectively. 

%\biggeoff{In the proposition below you should name the functors when you say ``extend to functors'' so that we know that you're defining presheaves of $\mathcal{A}$-objects $\mathcal{Z}F_n, \mathcal{B}F_n \in [(\tau^{op})^{op},\mathbf{Ch}(\mathcal{A})]$ so that it's clear that you really mean presheaves of chain complexes.}
\begin{prop}\label{cycle subobject functors as presheaves}
The cycle and boundary functors, $\cP \to \Sub_\cA(\cZ F_n^\fty)$, extend to presheaves of $\cA$-objects, $\tau^{op} \to \Sub_\cA(\cZ F_n^\fty)$, which we abusively also denote $\cZ F_n$ and $\cB F_n$. They are defined respectively by sending an arrow \begin{tikzcd}
	U \rar[tail, "\supset"] & V
\end{tikzcd} in $\tau^{op}$ to arrows, 

\begin{center}
\begin{tikzcd}
	\cZ F_n(U) \ar[dr, tail] \rar[rr, tail] &&  \cZ F_n(V) \ar[dl, tail] \\
	& \cZ F_n^\fty 
\end{tikzcd} \quad and \quad 
\begin{tikzcd}
	\cB F_n(U) \ar[dr, tail] \rar[rr, tail] &&  \cB F_n(V) \ar[dl, tail] \\
	& \cZ F_n^\fty 
\end{tikzcd}	
\end{center}
\end{prop}
\begin{proof}
We have shown each of the horizontal maps exist.
Identities and composition are preserved by the universal property of the limits involved.  
\end{proof}

\subsection{The Homological Memory Functor}\label{SS homological memory}

In this section we build a specific subobject to represent the $n$-cycles of $F$ which appear by a given up-set $U$ and which become ``filled-in'' as boundaries by a later up-set $V$.
We show this construction is functorial and call the induced functor $\mathcal{ZB}F_n(-,-)$.
It remembers the $n$'th homology features that have come by $U$ and gone by $V$, so we call it the \emph{homological memory functor}.

\begin{defn}\label{Def n-homological memory}
For $U \rightarrowtail V$ in $\tau^{op}$,\textit{the $n$'th homological memory} of the pair $(U,V)$ is the subobject $\ZB F_n(U,V) \rightarrowtail \cZ F_n^\fty$ induced by the following pullback square. 
\begin{center}
\begin{tikzcd}[column sep = huge, row sep = huge]
 \ZB F_n(U,V) \dar[tail] \rar[tail] \arrow[dr, phantom, "\usebox\pullback" , very near start, color=black] & \cB F_n(V) \dar[tail] \\
 \cZ F_n(U) \rar[tail] & \cZ F_n^\fty
\end{tikzcd}	
\end{center}
\end{defn}

\begin{prop}
$\ZB F_n : \mathbf{Arr}(\tau^{op}) \to \Sub_\cA( \cZ F_n^\fty)$ is a functor. 
\end{prop}
\begin{proof}
Since hom-sets in poset categories have at most one morphism, composition and identities are automatically preserved once we know how $\ZB F_n$ is defined on arrows $(U,V) \to (X,Y) \in \mathbf{Arr}(\tau^{op})$. 
Suppose $U \supseteq X \supseteq Y$ and $U \supseteq V \supseteq Y$ in $\mathbf{Arr}(\tau^{op})$. Then we have induced monics, 
\begin{center}
\begin{tikzcd}
	\cZ F_n(U) \rar[tail] & \cZ F_n(X) \end{tikzcd}\qquad \begin{tikzcd}
	\cB F_n(V) \rar[tail] & \cB F_n(Y)
\end{tikzcd},
\end{center}
making the bottom and right faces of the following cube commute: 
\begin{center}
\begin{tikzcd}
&\ZB F_n(X,Y) \arrow[dr, phantom, "\usebox\pullback" , very near start, color=black] \dar[dd, tail] \rar[rr, tail] && \cB F_n(Y) \dar[dd, tail] \\
\ZB F_n(U,V) \arrow[dr, phantom, "\usebox\pullback" , very near start, color=black] \ar[ur, dashed, tail]  \rar[rr, crossing over, tail] \dar[dd, tail] &&\cB F_n(V) \ar[ur, tail]  & \\
&\cZ F_n(X) \arrow[rr, tail] && \cZ F_n^\fty  \\
\cZ F_n(U) \ar[ur, tail] \rar[rr, tail] &&\cZ F_n^\fty \arrow[from =uu, crossing over, tail]  \ar[ur, equals] 
\end{tikzcd}. 
\end{center}
This makes $\ZB F_n(U,V)$ into (the vertex of) a cone of diagram that defines $\ZB F_n(X,Y)$ and induces the unique dashed arrow that makes the top and left faces of the cube commute by the universal property of the pullback, $\ZB F_n(X,Y)$. 
\end{proof}

\subsection{The Homological Lifespan Functor} \label{SS Lifespan subquot} 

For a filtration $F: \cP \to \mathbf{Ch}(\cA)$, the homological memory of the pair $U \supseteq V$ represents the $n$-cycles that have appeared by $U$ and that have become boundaries by $V$.

To extract the $n$-cycles that were born exactly at the initial boundary of $U$ and became boundaries exactly at the initial boundary of $V$, we want to discard all of the $n$-cycles that appeared by some open set before $U$ or that became boundaries at some open set before $V$. 

For this we use the blankets defined in Definition \ref{Def blanket}. Any $n$-cycles which were born before $U$ or became boundaries before $V$ must have been born by at least one of the blankets of $U$ or have become boundaries by at least one of the blankets of $V$ in $\tau^{op}$.
For readability we will refer to objects in $\mathbf{Arr}(\tau^{op})$ using capital letters $W,X,Y,Z$ and recall that they are in fact pairs of nested open sets in $\tau$ as necessary. 

\begin{defn}
For any object $X \in \mathbf{Arr}(\tau^{op})$, define $\Gamma_n(X)$ to be the cokernel in the following short exact sequence:

\begin{center}
\begin{tikzcd}
	0 \rar & \bigcup \limits_{W -\triangleright X}  \ZB F_n(W) \rar[] & \ZB F_n(X) \rar[] & \Gamma_n(X) \rar & 0
\end{tikzcd}.	
\end{center}
We also use the following quotient notation: 
\[ \Gamma_n (X)= \dfrac{\ZB F_n(X)}{\bigcup \limits_{W -\triangleright X}  \ZB F_n(W)}. \]
\end{defn}

\begin{rem} The image of the monomorphism on the left is equal to the kernel of the epimorphism on the right, and so the $n$-cycles represented in $\ZB F_n(W)$ for each \begin{tikzcd} W \rar[ -{Triangle[open]}] & X\end{tikzcd}, are being annihilated in $\Gamma_n(X)$ by the epimorphism on the right.
For $X = (U,V)$ where $U \supset V$ are open, we see that $\Gamma_n(X)$ represents the  $n$-cycles in $F$ whose 'lifespans' are precisely the set difference $U \setminus V \subseteq \cP$.
In this way, $\Gamma_n(X)$ represents the incremental change of $\ZB F_n$ at $X$, and so is reminiscent of a discrete derivative. 
\end{rem}

\begin{eg}
Let $K$ be a field and let $\cA = \mathbf{FinVect}_K$. 
Let $F: \mR \to \mathbf{FinVect}_K$ be a tame filtration.
After restricting to a finite sub-poset $\cP$ of $\mR$ and considering its topology of upper sets $\tau$, we recover finitely many half-open intervals $[i, j)$ for each pair of upper sets $[i, \fty) \supseteq [j, \fty)$ in $\mathbf{Arr}(\tau^{op})$.
Since every pair of elements in $\mR$ is comparable, there exist unique blankets $i' < i$ and $j' < j$ in $\cP$ provided that $i$ and $j$ are not initial in $\cP$.
Then 
\begin{align}\label{eqn Gamma[i,j)}
\Gamma_n([i,j)) = \dfrac{\cZ F_n([i,\fty)) \cap \cB F_n([j,\fty))}{(\cZ F_n([i',\fty)) \cap \cB F_n([j,\fty)) ) \cup (\cZ F_n([i,\fty)) \cap \cB F_n([j',\fty)) )}  
\end{align}
is the set of $n$-cycles which are born at $i$ and which bound an $(n+1)$-cycle by $j$.
This is precisely the pair-group of \cite[Section 2, Tame Functions]{E&H}. If this quotient is non-zero, then its rank counts the number of $n$-dimensional holes whose lifespan is the interval $[i,j)$.
By fixing a dimension, $n \in \mN$, and counting all of the lifespans of $n$-dimensional holes using the rank of $\Gamma_n([i,j))$ for each pair $[i, \fty) \supseteq [j, \fty)$ in $\mathbf{Arr}(\tau^{op})$, one recovers the data represented in a persistence diagram/barcode of the filtration $F$ as defined in \cite{CZ-ComputingMultiPersistence}. 
\end{eg}\

In the case of a single-parameter filtration, the rank of the quotient in Equation \ref{eqn Gamma[i,j)} can be computed using an inclusion-exclusion principle.
For filtrations indexed by two or more parameters, this is not the case, because the subobject lattices of objects in abelian categories are generally not distributive.
For example the subobject lattice of the Klein four group, $\mZ/2\mZ \oplus \mZ/2\mZ$, is a diamond lattice which is not distributive.
The focus of this paper is not to address this counting problem, but to elaborate on the idea that $\Gamma_n$ is a kind of functorial discrete derivative. 

\subsection{Change Action on Blankets}\label{SS blanket shift change action}

In this section show that every finite poset has an associated poset of, 'degree $n$ blankets,' for each $n \in \mN$, and see how $\mN^{op}$ acts on this poset by shifting 'degrees'.
This gives access to the elements that lie $n$-steps `below' a given element in a finite poset.
We use this in Section \ref{S CAD in multipersistence} to more precisely discuss the past homology cycles of a filtration for a given region, $V-U$. 

Let $\cP$ be a poset.
For each $x \in \cP$ define
\[ c^0(x) = \{ x \} \]
and inductively define  
\[ c^n(x) = \bigcup_{w \in c^{n-1}(x)} \{ y < w : y \leq z < w \implies y = z \}\]
to be the set of degree $n$ blankets of $x$.
To keep track of the root, $x$, define
\[ C^n(x) = c^n(x) \cup \{ x \}.\]
There is an induced ordering,
\[ C^n(x) \leq C^m(y) \text{ if and only if } x \leq y , n \geq m, \] 
which makes $\mathbf{P} = \{ C^n(x) : x \in \cP , n \in \mN \}$ a partially ordered set with an order preserving map to $\cP \times \mN^{op}$.

\begin{prop}\label{Prop inclusion into P times N}
There is a functor
\[  \begin{tikzcd} \mathbf{P} \rar["\iota"] & \cP \times \mN^{op} \end{tikzcd}\]
defined on objects by 
\[ \begin{tikzcd} C^n(x) \rar[mapsto] & (x,n)\end{tikzcd} .\]\
\end{prop}
\begin{proof}
The ordering on $\mathbf{P}$ coincides with the product partial order of $\cP \times \mN^{op}$, so $\iota$ is order preserving. 
\end{proof}
Note that when $c^n(x) = \varnothing$ we have $C^n(x) = \{ x \} = C^0(x)$ and $\iota(C^0(x)) = (x,0)$.
The monoid $\mN^{op}$, acts on $\mathbf{P}$ and $\cP \times \mN^{op}$, by addition in the second variable.
We use $\cP \times \mN^{op}$ to make formal definitions and proofs simpler but we think of the elements as blankets in $\mathbf{P}$. 

View $\mN$ as a poset category with opposite category $\mN^{op}$ and give $\mN \times \mN$ the product partial order.
Addition is an order preserving function $\mN \times \mN \to \mN$, so it's a functor between poset categories.
It's straightforward to check that $(\mN, +, 0)$ is a monoid in $\Cat$.
Since $(\mN \times \mN)^{op} = \mN^{op} \times \mN^{op}$, addition is order preserving on $\mN^{op}$ as well.
The proposition below follows. 

\begin{prop} \label{Prop N^op is a monoid} 
The triple, $(\mN^{op}, +, 0)$, is a monoid in $\Cat$. 	
\end{prop}
\begin{comment}
\begin{proof} 
 Addition,
 \[ + : \mN^{op} \times \mN^{op} \to \mN^{op} ; (n,m) \mapsto n+m,\]\
 
\noi  is an order preserving between partially ordered sets and therefore it's a functor between them when viewed as poset categories. The relevant monoid diagrams commute because ($\mN , + , 0$) is a monoid in $\Cat$.
 \end{proof}\
\end{comment}

\begin{prop}\label{Prop change action on blankets}
There is a change action, $( \cP \times \mN^{op} , \mN^{op}, 1_\cP \times (+), +, 0)$, in $\Cat$. 
\end{prop}
\begin{proof}
Let the monoid $(\mN^{op} , + , 0)$ act on $\cP \times \mN^{op}$ by addition in the second component. 
\begin{center}
\begin{tikzcd}[column sep = large]
\cP \times \mN^{op} \times \mN^{op} \ar[r, "1_\cP \times (+) "] & \cP \times \mN^{op} 
\end{tikzcd}.
\end{center}
\end{proof}

%Section 4 --------------------------------------------------------------------

\section{Change Action Derivatives in Multipersistence}\label{S CAD in multipersistence}

The rank of the pair group in \cite{E&H} counts the homology cycles whose lifespans are given by a certain interval.
Similarly, the rank of $\Gamma_n$ extracts information about the homology cycles whose lives span a region $V_U$ determined by a pair of nested opens $U \subseteq V$.
In this section we extend the rank of $\ZB F_n$ to a functor $ \mathbf{Arr}(\tau^{op}) \times \mN^{op} \to \mathbf{Arr}(G)$ for an abelian group $G$ and find the rank of $\Gamma_n$ lurking in a change-action derivative. 

Let $\cA$ be an abelian category, $\mathsf{K}(\mathcal{A})$ its Grothendieck group and let $\mathfrak{r} : \cA_0 \to \mathsf{K}(\mathcal{A})$ be the universal additive function that every other additive function into an abelian group factors through.
Let $\mathbf{Arr}(\mathsf{K}(\mathcal{A}))$ denote the arrow category of the delooping of $\mathsf{K}(\mathcal{A})$. 

\begin{prop}
The rank function induces a family of functors, 
\[ \begin{tikzcd}
\Sub_\cA(X)  \rar["\mathfrak{rk}_X "] & \mathbf{Arr}(\mathsf{K}(\mathcal{A})) \end{tikzcd}
 \] 
\noi defined by 
\[\left( \begin{tikzcd}
A \ar[dr, tail] \rar[rr, tail] && B \ar[dl, tail]  \\
\ & X &
\end{tikzcd} \right) 
\begin{tikzcd} \  \rar[mapsto] & \ \end{tikzcd} 
\left ( 
\begin{tikzcd}
* \dar["\mathfrak{r} (A) "'] \rar["0"] & * \dar[" \mathfrak{r} (B) "] \\
* \rar[" \mathfrak{r} (B/A) "'] & * 
\end{tikzcd} \right) \]

\end{prop}
\begin{proof} For any composable monomorphisms 
\[ \begin{tikzcd} A \rar[tail] & B \rar[tail] & C\end{tikzcd} \]
we have that 
\[ \mathfrak{r} (C/A) = \mathfrak{r}(C) - \mathfrak{r}(A) = \mathfrak{r}(C) -  \mathfrak{r}(B) + \mathfrak{r}(B)- \mathfrak{r}(A) =  \mathfrak{r}(C/B) + \mathfrak{r}(B/A) \]
showing the bottom of the diagram below commutes.  
\begin{center}
 \begin{tikzcd}[column sep = huge , row sep = large ]
                        * \dar["\mathfrak{r} (A) "'] \rar["0"] \rar[rr, bend left, "0"] & * \dar[" \mathfrak{r} (B) "] \rar["0"] & * \dar["\mathfrak{r}(C)"] \\
                        * \rar[" \mathfrak{r} (B/A) "'] \rar[rr, bend right, "\mathfrak{r}(C/A)"']  & *\rar["\mathfrak{r} (C/B)"']  & * 
\end{tikzcd} 
\end{center}
Identities are preserved because $\mathfrak{r}(A/A) = 0.$
\end{proof}

By Proposition \ref{Prop G1 subtraction change action in Cat}{,} $\mathbf{Arr}(\mathsf{K}(\mathcal{A}))$ acts on itself by subtraction.
That is, there is a change action $(\mathbf{Arr}(\mathsf{K}(\mathcal{A})) , \mathbf{Arr}(\mathsf{K}(\mathcal{A})), -, +, 0)$.
Recall that $c^n(U,V)$ denotes the $n$-th degree blankets of $(U,V)$ and let $\cP = \mathbf{Arr}(\tau^{op})$ as in the notation from the previous section.

\begin{prop}\label{Prop homological memory Extended}
For each $d \in \mN$ there is a functor, 
\[ \cup \ZB F_d : \mathbf{Arr}(\tau^{op}) \times \mN^{op} \to \Sub_\cA( \cZ F_d^\fty )\] 
defined by sending a morphism 
\[\begin{tikzcd}
((U,V), n) \rar[rrr, "\big( U \supseteq X{,} V \supseteq Y {,} n \geq m \big)  "] 
&&& ((X,Y), m)
\end{tikzcd}\]
to an inclusion of subobjects, 
\[
\begin{tikzcd}
 \bigcup \limits_{(W,Z) \in c^n(U,V)} \ZB F_d(W,Z)
\rar[ rr, tail, ""] && 
\bigcup \limits_{(W,Z) \in c^m(X,Y)} \ZB F_d(W,Z) 
\end{tikzcd}.
\]
\end{prop}
\begin{proof}
It is straightforward to see that identities are preserved and to see composition is preserved it suffices to show that for any degree $n$ blanket, $(W, Z) < (U,V)$, there exists a degree $m$ blanket $(W',Z') < (X,Y)$ with $(W, Z) \leq (W', Z')$.
The rest of the proposition follows from functoriality and the universal property of unions.\footnote{Note that the order here is the product partial order from $\tau^{op} \times \tau^{op}$, which is restriction of open subsets in each component.} 

%\biggeoff{I do like this proof, but the basic structure can be made a little clearer. I suggest saying that your first induction is a ``diagonal'' induction, i.e., you prove that given a degree $n$ blanket $(W,Z) < (U,V)$ there exists a degree $n$ blanket $(A,B) < (X,Y)$ with $(W,Z) < (A,B)$ (I know you used primes, but I'm being lazy in my TeX.). You can then indicate after establishing the diagonal case you induct on the $m$ component. The general math seems to parse, but explaining your technique will make it much clearer that the math, well, maths.}
We proceed by performing two inductions.
The first is an induction on the ``diagonal,'' that is we prove that given a degree $n$ blanket $(W,Z) < (U,V)$ there exists a degree $n$ blanket $(W',Z') < (X,Y)$ with $(W,Z) \leq (W',Z')$.
The base case, $n = 0,$ is trivial by taking $W = W'$ and $Z = Z'$.
Now assume that the result holds for $n = m = k$ for some $k > 0$, and consider $n = m = k + 1$.
As we saw in Remark \ref{rem characterizing blankets in Arr(tauo^p)}, every degree $1$ blanket of $(U,V)$ is either of the form $(W, V)$, for some blanket $W < U$, or $(U, Z)$, for some blanket $Z < V$; similarly for $(X,Y)$.
For a degree $1$ blanket, $(W, V) < (U, V)$, we have the following cases:

\begin{itemize}
    \item The case when $U = X$ and $V = Y$ is trivial.
    
    \item If $U = X$ and $V \neq Y$, then the blanket $(W,Z) < (U,V)$ must have that $W < U = X$ is a blanket. Similarly there exists a blanket $Z < Y$ such that $V \leq Z < Y$. This means $(X, Z) < (X,Y)$ is a blanket and we have 
    \[ (W, V) < (U, V) = (X,V) \leq (X, Z) < (X,Y) \]
    Applying the induction hypothesis for degree $k$ blankets to each of gives the result for a subset of the degree $k + 1$ blankets. 
    
    \item If $U \neq X$ and $V = Y$ then, there exists a blanket $U' < X$ with with $U \leq U'$. Since the blanket $(W,V) < (U,V)$ we have 
    \[ (W,V) < (U, V) \leq (U', V) = (U', Y) \leq X,Y) \]
    Applying the induction hypothesis for degree $k$ blankets to each of these gives the result for a subset of the degree $k + 1$ blankets. 
\end{itemize}

The argument for blankets of the form $(U, Z) < (U, V)$ is similar.
Every degree $k + 1$ blanket is a degree $k$ blanket of a degree $1$ blanket, and the cases above together account for the degree $k$ blankets of all different possible degree $1$ blankets.
This concludes the first induction: the result holds for all $n = m \in \mN$. 

For the second induction we fix $m \in \mN$ and proceed by induction on $k$, where $n = m + k$.
The base case, $k = 0$, follows from above, so we assume the result holds for $n = m+k$ for some $k > 0$.
For $n = m + k + 1$, notice that each degree $1$ blanket of $(U,V)$ has the form $(W, V)$ or $(U, Z)$, we can apply the induction hypothesis as we did in the previous induction. 

Both of these inductions together show that for every $n \geq m$ in $\mN$ and $(U, V) \leq (X,Y)$, for every degree $n$ blanket $(W, Z) < (U, V)$ there exists a degree $m$ blanket $(W', Z') < (X,Y)$, and a restriction $(W, Z) \leq (W', Z')$.
As mentioned at the start, the proof of the proposition follows from functoriality of $\ZB F_d$ and the universal property of unions. 

%This proof isn't quite right. 
\begin{comment}
Assume $U \supseteq X , V \supseteq Y $, and $n \geq m$. 

%This part here isn't necessarily true.
Then $U \in C^k(X)$ and $V \in C^\ell(Y)$ for some $k, \ell \in \mN$. Assume $k \geq  \ell a$, without loss of generality. 

For any $(W,Z) \in c^n(U,V)$ we then have that $W \in C^{n + k(X)}$ and $Z \in C^{n + \ell}(Y)$. 

Since $n+k \geq n + \ell \geq n \geq m$ so we get an inclusion 
\[ \begin{tikzcd}
\ZB F_d(W,Z) \rar[tail] & \ZB F_d( W' , Z')\end{tikzcd}\]

\noi of subobjects in $\Sub_\cA( \cZ F_d^\fty)$ for some $(W',Z') \in c^m(X,Y)$. Doing this for all $(W,Z) \in c^n(U,V)$ induces the desired subobject inclusion by the universal property of unions. 
\end{comment} 
\end{proof}\

\begin{defn}
The \textit{rank of $\cup \ZB F_d$} is denoted $\zb F_d$ and is defined by post-composing with the rank functor $\mathfrak{rk}$. 

\[\begin{tikzcd}[column sep = large, row sep = large] 
	    	\mathbf{Arr}(\tau^{op}) \times \mN^{op}  \rar["\cup \ZB F_d"] \ar[dr, "\mathbf{zb} F_d"']  & \Sub_\cA(\cZ F_d^\fty) \dar["\mathfrak{rk}"] \\
	       	\ & \mathbf{Arr}(\mathsf{K}(\mathcal{A})) 
            \end{tikzcd}\]
\end{defn}

The following theorem is the main result of this paper which allows us to extract a change-action derivative from the functor $\cup \ZB F_d$. 

\begin{thm}\label{thm extended homological memory functor is differentiable} $\mathbf{zb} F_d$ has a change-action derivative in $\Cat$
 
\[ \begin{tikzcd} \mathbf{Arr}(\tau^{op}) \times \mN^{op} \times \mN^{op} \rar["\mathbf{\Gamma}_d "] & \mathbf{Arr}\left( \mathsf{K}(\mathcal{A}) \right)\end{tikzcd} .\]
\end{thm}
\begin{proof}
As every functor from a category with a change action structure into the arrow category of an abelian group is change-action differentiable, this follows immediately from Proposition \ref{thm functors into G1 are differentiable}. 
\end{proof}

%\biggeoff{Two main notes. First, you need to kill the \texttt{\textbackslash newpage} so that the reader can see the paper continues. On a style note, I find that it's not great form to do page breaks unless it's between sections, subsections, or chapters. I've commented the \texttt{\textbackslash newpage} out, but if you disagree with me you should kill this comment and the \% tag.

%Second, it would be helpful for your readers to indicate why it is that this theorem and your next corollary are so cool! In particular, you just showed that these persistence homological gadgets are recovered by your categorical techniques AND by change action derivative stuff, and emphasizing this to your readers is helpful because they may not catch this!}

As a direct consequence we obtain the following corollary, which shows how the change-action derivative determined by categorical calculus of finite difference, Theorem \ref{thm extended homological memory functor is differentiable}, captures a functorial generalization of the rank of the pair-group of \cite{E&H}, $\mathfrak{r} \circ \Gamma_d$. 

\begin{cor}\label{cor lifespan functor recovered by derivative}
The functor $\mathfrak{r}\circ \Gamma_d$ can be recovered by precomposing the change-action derivative $\mathbf{\Gamma}_d$ with the embedding 
\[ \begin{tikzcd} \mathbf{Arr}(\tau^{op}) \rar[tail, "i_{0,1}" ] & \mathbf{Arr}(\tau^{op}) \times \mN^{op} \times  \mN^{op} \end{tikzcd}.\]
defined by
\[\begin{tikzcd} (U,V) \rar[mapsto] &  \big((U,V) , 0, 1 \big) \end{tikzcd}  \]
\end{cor}
\begin{proof}
For any $(U,V)  \in \tau^{op}_1$, 
\begin{align*} 
i_{0,1} \mathbf{\Gamma}_d (U,V) 
= \mathbf{\Gamma}_d \big( (U,V) , 0, 1 \big) 
&= \mathbf{zb}F_d((U,V),0)  - \mathbf{zb}F_d((U,V),1) \\
&= \mathfrak{r} \left( \ZB F_d( U,V ) \right) - \mathfrak{r} \left( \bigcup_{(W,Z) \in c^1(U,V)} \ZB F_d ( W,Z) \right)  \\
&= \mathfrak{r} \left( \dfrac{ \ZB F_d( U,V ) }{ \bigcup \limits_{(W,Z) \in c^1(U,V)} \ZB F_d ( W,Z) } \right) \\
&= \mathfrak{r} \Gamma_d ( U,V )
\end{align*}
where the third equality comes from the definition of the rank function.
\end{proof}

%Random Thoughts
\begin{comment} 
- there's more data in this change-action derivative than just the data of the persistence diagram; can it be reorganized in terms of a functor-calculus kind of approximation? 
    - functor calculus is a homotopy theoretic thing, and the issue with a complete invariant for multipersistence could be that the notion of sameness is too strict there. Kristine's crew has a characterization of homotopy colimits in terms of functor calculus that involves monoidal axtions in Cat and I'm exploring the idea of the persistence diagram really being a kind of homotopy colimit. 
    -Main Idea here: no complete invariants for filtrations in $\mR^n$ for $n > 1$ exist...is the issue homotopy theoretic because we have room to wiggle lifespans around? If it is, then can we expect a complete homotopical invariant? This likely won't be `stable under small perturbations,' but can it still be useful? 
\end{comment} 

%\nocite{*}
\bibliographystyle{amsplain}
\bibliography{biblio}

@incollection {E&H,
    AUTHOR = {Edelsbrunner, H. and Harer, J.},
     TITLE = {{P}ersistent {H}omology---{A} {S}urvey},
 BOOKTITLE = {Surveys on discrete and computational geometry},
    SERIES = {Contemp. Math.},
    VOLUME = {453},
     PAGES = {257--282},
 PUBLISHER = {Amer. Math. Soc., Providence, RI},
      YEAR = {2008},
   MRCLASS = {55N35 (52-02 55T05 57M99 68W01)},
  MRNUMBER = {2405684},
MRREVIEWER = {Greg Friedman},
       DOI = {10.1090/conm/453/08802},
       URL = {https://doi.org/10.1090/conm/453/08802},
}

@article{Cai2013ATO,
  title={A theory of changes for higher-order languages: incrementalizing $\lambda$-calculi by static differentiation},
  author={Yufei Cai and Paolo G. Giarrusso and Tillmann Rendel and Klaus Ostermann},
  journal={Proceedings of the 35th ACM SIGPLAN Conference on Programming Language Design and Implementation},
  year={2013},
  url={https://api.semanticscholar.org/CorpusID:196189637}
}

@book {CFTWM,
    AUTHOR = {Mac Lane, Saunders},
     TITLE = {Categories for the working mathematician},
    SERIES = {Graduate Texts in Mathematics},
    VOLUME = {5},
   EDITION = {Second},
 PUBLISHER = {Springer-Verlag, New York},
      YEAR = {1998},
     PAGES = {xii+314},
      ISBN = {0-387-98403-8},
   MRCLASS = {18-02},
  MRNUMBER = {1712872},
}

@article{AP&L,
   title={Cartesian Difference Categories},
   volume={Volume 17, Issue 3},
   ISSN={1860-5974},
   url={http://dx.doi.org/10.46298/lmcs-17(3:23)2021},
   DOI={10.46298/lmcs-17(3:23)2021},
   journal={Logical Methods in Computer Science},
   publisher={Centre pour la Communication Scientifique Directe (CCSD)},
   author={Alvarez-Picallo, Mario and Lemay, Jean-Simon Pacaud},
   year={2021},
   month=sep }

@article {CZ-ComputingMultiPersistence,
    AUTHOR = {Carlsson, Gunnar and Singh, Gurjeet and Zomorodian, Afra},
     TITLE = {Computing multidimensional persistence},
   JOURNAL = {J. Comput. Geom.},
  FJOURNAL = {Journal of Computational Geometry},
    VOLUME = {1},
      YEAR = {2010},
    NUMBER = {1},
     PAGES = {72--100},
   MRCLASS = {52B55 (68Q17)},
  MRNUMBER = {2770959},
       DOI = {10.20382/jocg.v1i1a6},
       URL = {https://doi.org/10.20382/jocg.v1i1a6},
}

@inproceedings{peristenceBarcodesforShapes, author = {Carlsson, Gunnar and Zomorodian, Afra and Collins, Anne and Guibas, Leonidas}, title = {Persistence barcodes for shapes}, year = {2004}, 
isbn = {3905673134}, 
publisher = {Association for Computing Machinery}, 
address = {New York, NY, USA}, 
url = {https://doi.org/10.1145/1057432.1057449},
doi = {10.1145/1057432.1057449}, 
booktitle = {Proceedings of the 2004 Eurographics/ACM SIGGRAPH Symposium on Geometry Processing}, 
pages = {124–135}, 
numpages = {12}, 
location = {Nice, France}, 
series = {SGP '04} }

\end{document}